\documentclass{amsart}
\usepackage{amsmath}
\usepackage{amsthm}
\usepackage{amssymb}
\usepackage{amsfonts}
\usepackage{amsxtra}     
\usepackage{epsfig}
\usepackage{graphicx}
\usepackage{graphics}
\usepackage[table]{xcolor} 
\usepackage{longtable}
\usepackage{multicol}
\usepackage{hyperref}
\usepackage{array}
\newsavebox\ltmcbox
\newcounter{entryno}
\def\tabline{\the\value{entryno} & Description\addtocounter{entryno}{1}\\}

\newtheorem{theorem}{Theorem}

\DeclareMathOperator{\vol}{vol}
\newtheorem{lemma}[theorem]{Lemma}
\newtheorem{corollary}[theorem]{Corollary}
\newtheorem{sub-lemma}[theorem]{Sub-Lemma}

\author{Taiyo Inoue}
\title{The 825 Smallest Right-Angled Hyperbolic Polyhedra}
\email{taiyoinoue@gmail.com}

\begin{document}

\maketitle

\begin{abstract} An algorithm for determining the list of smallest volume right-angled hyperbolic polyhedra in dimension 3 is described.  This algorithm has been implemented on computer using the program {\tt Orb} to compute volumes, and the first 825 polyhedra in the list have been determined.  
\end{abstract}

\section{Introduction}
In a prior paper \cite{I}, I described how to organize the volumes of the infinite family of compact, right-angled hyperbolic polyhedra in hyperbolic $3$-space $\mathbb{H}^3$.  This involves applying one of two combinatorial operations to the 1--skeleton of the polyhedron -- either delete an edge (edge-surgery) or split the polyhedron along an ``incompressible'' polygon into a pair of polyhedra (decomposition).  The effect of this modification is to reduce the average complexity of the polyhedron, while keeping it in the family of right-angled polyhedra.  Repeated application of these modifications produces a chain of polyhedra which terminates in a family of right-angled polyhedra which are ``atomic'' with respect to these operations.  These are called L\"obell polyhedra.  Every non-L\"obell right-angled hyperbolic polyhedron can have one of these operations performed on it.  In particular, every right-angled hyperbolic polyhedron can be obtained by applying these operations in reverse -- either add an edge or glue two polyhedra together.  

The geometry of a right-angled hyperbolic polyhedron is determined completely by its combinatorics.  Once the 1--skeleton of the polyhedron is given (a trivalent planar graph with some additional conditions to be described), the geometric structure which makes the polyhedron right-angled and hyperbolic is unique in the sense that any other right-angled hyperbolic polyhedron with an isomorphic 1--skeleton is isometric to the one given.  This leads one to suspect that a combinatorial approach to geometric problems, such as the ordering of volumes, is feasible.

In fact, the combinatorial operations of edge-deletion and decomposition described above have realizations as geometric operations.  Edge-deletion looks like increasing the dihedral angle measure of the edge being deleted so that it goes from $\frac{\pi}{2}$ to $\pi$ while also leaving all the other edges with right-angled dihedral angles -- so a sort of continuous flattening.   This deformation is analogous to 3--manifold surgery.  Decomposition is like cutting all the way through the polyhedron with a knife, then molding of the resulting pieces to make them right-angled polyhedra.

It happens that each of these operations decreases the volume of the polyhedron, except in some special cases of decomposition where the splitting surface is a totally geodesic polygon (in this case, the volume does not change).  It is also the case that the average number of faces of each component of the polyhedron will decrease (keep in mind that decomposition actually splits a polyhedron into pieces and I am blurring the distinction between a polyhedron (singular) and polyhedra (plural)).  Therefore, this process will eventually terminate.  In fact, this process terminates in a family of L\"obell polyhedra.

So we have the following situation.  Start with a non-L\"obell right-angled hyperbolic polyhedron $P$.  Repeatedly apply either edge surgery or decomposition resulting in a chain of right-angled polyhedra eventually terminating in a collection of L\"obell polyhedra $L$:

$$P \rightarrow P_1 \rightarrow P_2 \rightarrow ... \rightarrow P_k = L$$ 

Taking the volume of each polyhedron in the above chain produces a chain of inequalities of volumes of right-angled hyperbolic polyhedra:

$$\text{Vol}(P) \geq \text{Vol}(P_1) \geq \text{Vol}(P_2) ... \geq \text{Vol}(P_k) = \text{Vol}(L)$$ 
The number $\text{Vol}(L)$ can be explicitly calculated since there are known formulas for the volumes of right-angled hyperbolic L\"obell polyhedra \cite{V}.

A simple corollary of this result is that the smallest and second-smallest right-angled hyperbolic polyhedra are the L\"obell polyhedra $L_5$ and $L_6$ respectively.  The 1--skeleta of these polyhedra are in Figure 1.  Note that $L_5$ is the right-angled dodecahedron.  The main result described in this paper is to extend this list from the two smallest polyhedra to the 825 smallest polyhedra via an algorithm implemented on computer.  Also provided is a way to extend this list to many thousands of entries if certain computational obstacles can be overcome.

\begin{figure}
\centering
\includegraphics[height=4cm]{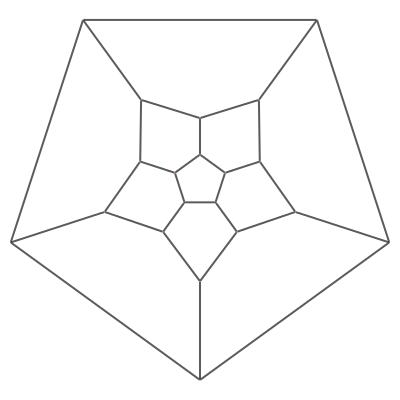}\includegraphics[height=4cm]{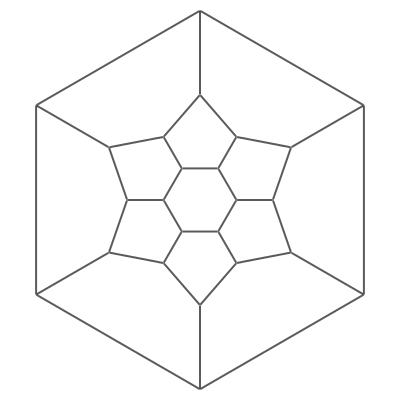}\includegraphics[height=4cm]{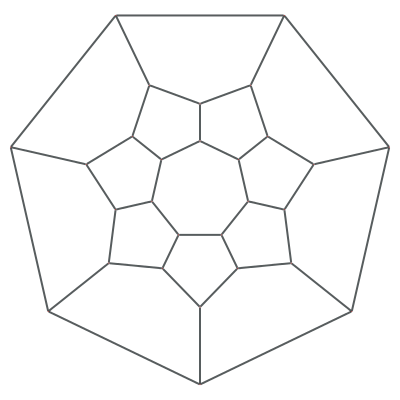}
\includegraphics[height=4cm]{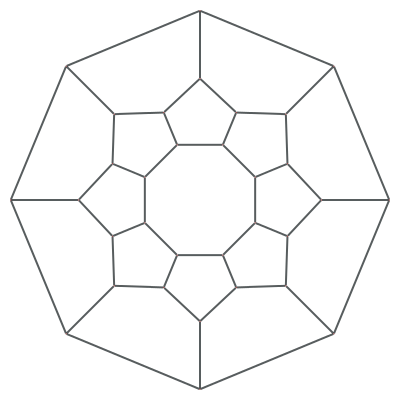}\includegraphics[height=4cm]{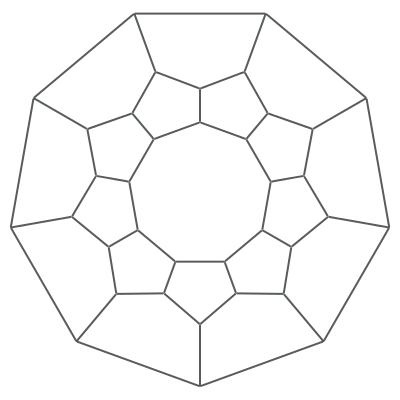} \includegraphics[height=4cm]{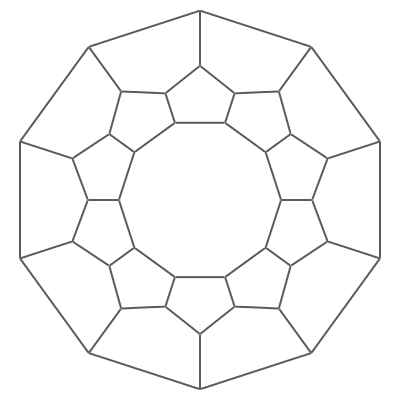}
\caption{Planar projections (Tutte embeddings) of the 1--skeleta of the first six L\"obell polyhedra.}
\label{lobells}
\end{figure}

\section{Preliminaries}

The setting for the polyhedra throughout this paper is hyperbolic 3--space ${\mathbb{H}}^3$.  A \emph{right-angled hyperbolic polyhedron} is a compact intersection of finitely many hyperbolic half-spaces such that if the boundaries of two of these half-spaces intersect, then the intersection has a dihedral angle measure of $\frac{\pi}{2}$.  Such polyhedra can be viewed as Coxeter orbifolds. The action of the group generated by reflections in the boundaries of half-spaces of right-angled polyhedra produces beautiful tilings of $\mathbb{H}^3$ similar to the familiar tiling of euclidean space by cubes. \footnote{See, for example, the video \emph{Not Knot} for a guided tour of $\mathbb{H}^3$ tiled by right-angled hyperbolic dodecahedra \cite{NotKnot}.}

In 1967, A. Pogorelov classified right-angled hyperbolic polyhedra by the combinatorics of their 1--skeleta \cite{Pog}:

\smallskip

\noindent \textbf{Theorem:} (Pogorelov 1967) \emph{A polyhedron $P$ has a geometric realization in $\mathbb{H}^3$ as a right-angled hyperbolic polyhedron if and only if:
\begin{enumerate}
 \item The 1--skeleton of $P$ is trivalent.
 \item There are no prismatic 3 or 4--circuits (defined below).
\end{enumerate}
This geometric realization is unique up to isometry.} 

\smallskip
A \emph{prismatic $k$--circuit} is a closed curve on the boundary of the polyhedron which intersects exactly $k$ edges transversely such that no two of these edges share a vertex.  A prismatic $k$--circuit can be viewed as a polyhedral analog of the notion of incompressible surfaces for 3--manifolds as it bounds a topological $k$--gon in the polyhedron.  If the hyperbolic polyhedron is right-angled, we can understand the prohibition of prismatic 3 and 4--circuits as being analogous to the prohibition of essential spheres and irreducible tori in hyperbolic 3--manifolds since right-angled triangles are naturally spherical and right-angled quadrilaterals are euclidean.

Pogorelov's result demonstrates that right-angled hyperbolic polyhedra are fundamentally combinatorial as they are determined completely by their 1--skeleton.

Two distinct faces $F_1$ and $F_2$ of a right-angled hyperbolic polyhedron $P$ are \emph{edge-connected} if they are nonadjacent and there exists an edge $e$ of $P$ connecting a vertex of $F_1$ to a vertex of $F_2$.  Such an edge will be said to \emph{edge-connect} $F_1$ and $F_2$.  Since $P$ has a trivalent 1--skeleton, every edge of $P$ edge-connects a unique pair of faces.  

A face $F$ of $P$ is \emph{large} if it is a $k$--gon with $k \geq 6$.  If an edge $e$ edge-connects two large faces and it is not intersected by a prismatic 5--circuit, then this edge is said to be \emph{very good}.  It has been shown that deleting the interior of a very good edge and demoting the endpoints to non-vertices results in a new right-angled hyperbolic polyhedron with one less face, three fewer edges, and two fewer vertices. This new polyhedron has strictly less volume than the one started with.  Call this process \emph{edge-deletion}.

If a polyhedron $P$ has a large face $F$, then we can add an edge to $P$ by connecting two points which lie on the interiors of two different edges, as long as these edges are separated by at least two edges in $F$.  These conditions ensure that the polyhedron which results will not have any prismatic 3 or 4--circuits, and so will have a realization as a right-angled hyperbolic polyhedron.   Call this operation \emph{edge-addition}.  

Every right-angled hyperbolic polyhedron other than the dodecahedron has a large face.  Therefore, edges can almost always be added to a polyhedron to generate new examples of right-angled polyhedra.  Call the polyhedra which result from edge-addition \emph{edge-children} and the set of all edge-children, edge-children of the edge-children, etc. the \emph{edge-descendants}.

The other fundamental operation needed is \emph{composition}, and its inverse, \emph{decomposition}.  If two polyhedra $P_1$ and $P_2$ both have a face which is a $k$--gon, then by choosing a combinatorial isomorphism of these $k$--gons, $P_1$ and $P_2$ can be glued together along these faces.  Most of the time, the faces will not be isometric to one another in the geometric realization of $P_1$ and $P_2$ and so this gluing may not make so much sense geometrically.  But it's certainly fine to do topologically or combinatorially -- identify the subgraphs of the 1--skeleta corresponding to the $k$--gons using the combinatorial isomorphism, delete the interiors of the edges of the $k$--gons, and demote the vertices to non-vertices.  It turns out that this operation, which will be called \emph{composition}, results in a right-angled hyperbolic polyhedron. The inverse operation, splitting a polyhedron into two right-angled polyhedra, is called \emph{decomposition}.

Decomposition should be viewed as the polyhedral analog of decomposition of Haken hyperbolic 3--manifolds along incompressible surfaces.  Note that composition along a $k$--gon will have a distinguished prismatic $k$--circuit bounding a topological $k$--gon in the polyhedron which is analogous to an incompressible surface. Taking this analogy a little further, a result of Agol, Storm, Thurston and Dunfield \cite{ADST} sheds light on the effect of composition/decomposition on volume. If $P$ is a composition of $Q_1$ and $Q_2$, then 
$$\vol(P) \geq \vol(Q_1) + \vol(Q_2)$$ 
with equality occurring if and only if the polygonal gluing sites are isometric.

An infinite family of right-angled hyperbolic polyhedra called the \emph{L\"obell polyhedra} can be constructed as follows.  Start with an $n$--gon where $n \geq 5$.  Attach to each edge of this $n$--gon a pentagon, and glue the remaining edges of the pentagon together to produce a bowl shape.  Two copies of this can be glued together to produce a polyhedron $L_n$ with two $n$--gons, and $2n$ pentagons. These polyhedra are notable in that they have no very good edges and are not formed by the composition of any two right-angled polyhedra.
See Figure \ref{lobells} for pictures of the L\"obell polyhedra.

The main result of \cite{I} is the following: 

\begin{theorem} {\sl Let $P_0$ be a compact right-angled hyperbolic polyhedron.  Then there exists a sequence of disjoint unions of right-angled hyperbolic polyhedra $P_1, P_2, \dots, P_k$ such that for $i=1,\dots, k$, $P_i$ is gotten from $P_{i-1}$ by either a decomposition or edge surgery, and $P_k$ is a set of L\"obell polyhedra. Furthermore, 
$$\vol(P_0) \geq \vol(P_1) \geq \vol(P_2) \geq \dots \geq \vol(P_k).$$}
\end{theorem}
\begin{corollary} Every right-angled hyperbolic polyhedron can be obtained by a sequence of either edge additions or compositions applied to a set of L\"obell polyhedra.
\end{corollary}

\section{Description of the Algorithm}

{\tt Orb} is a program written by Damien Heard for computing geometric structures of many 3--dimensional orbifolds \cite{Heard}.  It is an extension of Jeff Weeks' {\tt Snappea} program for 3--manifolds \cite{weeks2001snappea}.  Orbifolds with underlying space $S^3$ are typically entered into {\tt Orb} by manually drawing the singular locus of the orbifold on a planar canvas that is part of the user interface.  The singular locus will be a (possibly) knotted graph.  The user then specifies the cone angle at each edge.  {\tt Orb} triangulates the resulting orbifold and numerically computes the hyperbolic structure on each tetrahedron which glue together to produce the orbifold (assuming a hyperbolic structure exists).  One can then probe various aspects of the geometrized orbifold such as volume, lengths of geodesics and many other things.  

{\tt Orb} can compute the volume of a right-angled hyperbolic polyhedron $P$ by constructing an orbifold double cover $Q$.  This orbifold is easily visualized by taking two copies of $P$ and gluing each face of one copy to the corresponding face of the other.  As $P$ is homeomorphic to a ball, this doubling operation produces an orbifold $Q$ with underlying space $S^3$ whose singular locus is the 1--skeleton of $P$. The cone angle of each edge of the singular locus of $Q$ is $\pi$.  This orbifold $Q$ has a hyperbolic structure, and since it is a double cover of $P$, the volume of $Q$ is precisely twice that of $P$.  
All volumes reported in this paper come from {\tt Orb} in this way.

The orbifolds constructed in {\tt Orb}'s user interface can be saved to a plain text file which stores the triangulation of the orbifold in Casson format.  If the geometric structure of the orbifold has been computed, then this data is also saved to the file.

The algorithm used here implements a highly stripped-down version of {\tt Orb} called {\tt Vol} (named simply to distinguish from {\tt Orb}) which works from a command line interface.  {\tt Vol} accepts an orb file and outputs the volume to the console.  This simplified console version of {\tt Orb} the makes automating the computation and ordering of volumes easier.

Edge-addition is implemented as a {\tt python} program.  It accepts an {\tt Orb} file containing a triangulation of $P$ as input.  The program will then find large faces of $P$ and add an edge connecting points of two edges which separated by at least two edges on either side.  It will then triangulate the resulting orbifold to produce an edge-child of $P$, then save the result as an {\tt Orb} file.  This function can be repeatedly invoked to produce all of the edge-children of $P$.

I did not bother to implement composition of polyhedra as a computer program.  The justification for this is given below.

The algorithm proceeds as follows.  A file $\mathcal{L}$ is read at the start.  It serves as a plain text database containing a row for each polyhedron that has been constructed.  Each invocation of the algorithm can extend the number of polyhedra represented in $\mathcal{L}$ by quite a large number. However, each invocation extends the list of smallest volume polyhedra by exactly one. This is because $\mathcal{L}$ contains the rank-ordered list of smallest polyhedra, but also contains the children of these.

If the algorithm has never been invoked on $\mathcal{L}$, then $\mathcal{L}$ contains only certain polyhedra chosen to generate a complete list of the polyhedra whose volume is below some pre-determined value (more on this in Theorem \ref{fifteen}).  These initial polyhedra serve as seeds for the program to work from.  It is clear that some of the smaller-volume L\"obell polyhedra must be included in this initial list, but there are others that must be added as well.  More on this later.

\begin{figure}
\centering
 \includegraphics[trim={5.2cm 12cm 0 6cm}, clip]{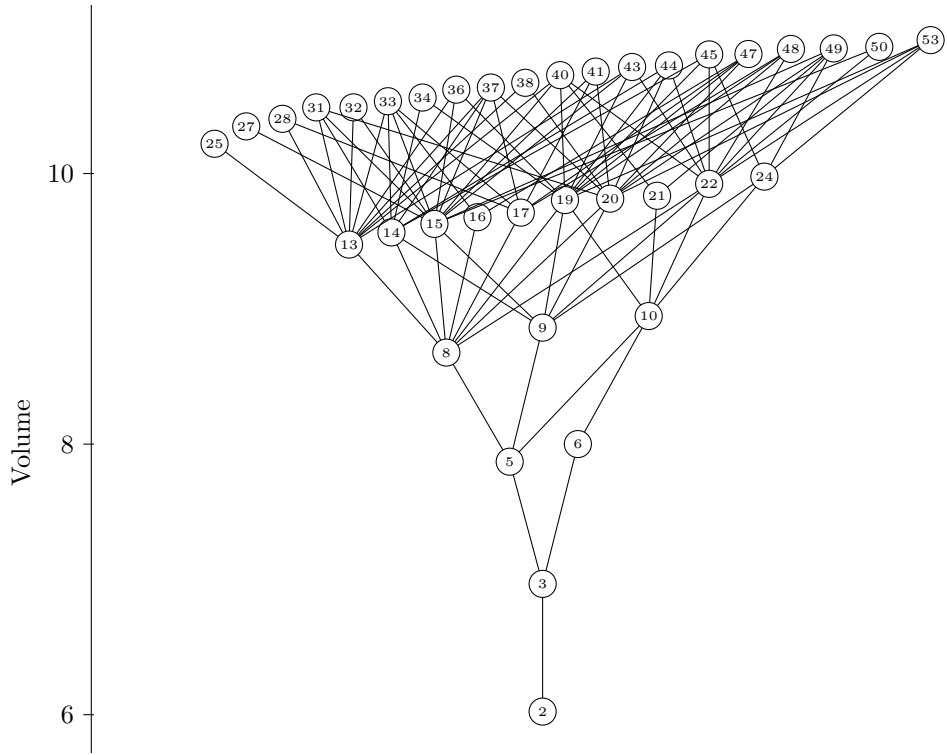}
 \caption{The first five generations of edge-descendants of $L_6$ (labeled ``2'').  The numerical labels indicate rank ordering in the list of smallest volume polyhedra.} 
 \label{familytree}
\end{figure}

The algorithm scans the list $\mathcal{L}$ by volume and finds the smallest polyhedron $P$ whose edge-children have not yet been computed.  This polyhedron is the next-smallest polyhedron in the list and is recorded as such by assigning it its rank-ordering in terms of smallest volume (for example, $L_5$ is identified as ``1.orb'' in the list $\mathcal{L}$ after the very first invocation of the algorithm).  When it finds $P$, the algorithm computes every possible child of it.  Each child $C$ is then tested to see if it already on the list or not.  This is done by checking if the 1--skeleton of $C$ (a trivalent planar graph) is graph-isomorphic to the 1--skeleton of a polyhedron already on the list.  

For the graph isomorphism task, the {\tt python} module {\tt{graph-tool}} is used \cite{Peixoto}.  It contains a function for testing graph isomorphism.  However, for the sake of efficiency, before this {\tt graph-tool} function is called, the \emph{face vector} of the child $C$ is computed.  This face vector is a vector of integers which keeps track of the number of pentagons, hexagons, septagons, octagons, etc. among the faces of $P$.  So, for instance, $L_8$ has face vector $[16, 0, 0, 2]$ since it has 16 pentagonal faces and 2 octagonal faces.  From $\mathcal{L}$, all of the polyhedra which have the same face vector as the child are collected -- obviously those which do not have the same face vector as $C$ cannot be isometric to $C$.  This serves as a computationally cheap first-pass filter which cuts down on the number of graph isomorphism checks that need to be done.

However, there can be hundreds of polyhedra listed in $\mathcal{L}$ which have the same face vector as the child $C$.  In this case, rather than calling the graph isomorphism function hundreds of times, the algorithm will instead run the child $C$ through {\tt Vol}, then collect from $\mathcal{L}$ all the polyhedra with the same recorded volume as $C$.  Then the graph isomorphism test compares $C$ with this presumably smaller family of polyhedra.  If $C$ is not found to be in the list of polyhedra with the same volume, then it is tested for graph isomorphism against the entire list of polyhedra with the same face vector.  This is to ensure that the absence of $C$ from the list of same-volume polyhedra isn't due to numerical error in {\tt Vol}.  

At this point, if $C$ is not found to be in the list $\mathcal{L}$, it is run through {\tt Vol} (if it has not already) to obtain its volume.  A new row is added to $\mathcal{L}$ containing the relevant information.  If $C$ is found in $\mathcal{L}$, then the fact that $P$ is a parent of $C$ is recorded in the row corresponding to $P$ so as to record the ``family tree'' of right-angled polyhedra.  See Figure \ref{familytree}.

But there is an additional wrinkle in this process.  Sometimes for particular triangulations of polyhedra with a relatively large number of faces (25 or greater), {\tt{Vol}} is unable to find a geometric solution to the gluing equations.  A rejiggering of the triangulation can sometimes get {\tt{Vol}} to find the geometric solution and hence the correct volume.  However, there are still a handful of polyhedra which do not have an ``{\tt{Vol}}-certified'' volume.  

This is an obstacle to getting a long list of smallest-volume polyhedra. Examining the parentage of the problematic polyhedra reveals that this obstacle will limit the correctness of the list to the first 825 right-angled polyhedra.  If this problem could be fixed, then the list of smallest-volume polyhedra could go from 825 polyhedra to many thousands of polyhedra long.

To summarize:

\renewcommand{\labelitemii}{$\bullet$}
\renewcommand{\labelitemiii}{$\bullet$}
\renewcommand{\labelitemiv}{$\bullet$}

\begin{enumerate}
 \item[(1)] A file $\mathcal{L}$ which contains a list of initial polyhedra (the seeds) and some of their edge-descendants is read into memory.  
 \item[(2)] The smallest volume polyhedron $P$ listed in $\mathcal{L}$ whose edge-children have not yet been computed is found.  This is the next smallest volume polyhedron.  The orb file which triangulates $P$ is read into memory.
 \item[(3)] A {\tt python} program finds the large faces of $P$ and performs all possible edge-additions until all edge-children have been created.
 \item[(4)] Each edge-child $Q$ is compared against the polyhedra in the list $\mathcal{L}$ to see whether it can be found in $\mathcal{L}$.
 \begin{itemize}
  \item To reduce the number of polyhedra in $\mathcal{L}$ which need to be compared to $Q$, polyhedra in $\mathcal{L}$ with the same face vector as $Q$ are collected.  Call this list $\mathcal{F}$.
      \begin{itemize}
	\item If the number of polyhedra in $\mathcal{F}$ is bigger than some threshold value, then the volume of $Q$ is computed using {\tt Vol} and all polyhedra in $\mathcal{L}$ with the same computed volume are collected.  Call this list $\mathcal{V}$.
	  \begin{itemize}
	    \item {\tt graph-tool} tests for graph isomorphism between $Q$ and elements of $\mathcal{V}$. 
	      \begin{itemize}
		\item If a graph isomorphism is found, go to step (5a).
		\item If a graph isomorphism is not found, then $Q$ is tested for graph isomorphism against all the elements in $\mathcal{F}$.  If it is not found to be isomorphic to something in $\mathcal{F}$, go to step (5b).  Otherwise, go to step (5a).
	      \end{itemize}
	  \end{itemize}
	 \item If the number of polyhedra in $\mathcal{F}$ is smaller than some threshold value, then each polyhedron in $\mathcal{F}$ is tested to see if it is isomorphic to $Q$. 
	\end{itemize}
 \end{itemize}
  \item[(5a)] If $Q$ is found to be in $\mathcal{L}$, then the fact that $P$ is a parent of $Q$ is recorded in the row corresponding to $P$, and $Q$ is then discarded.
 \item[(5b)] If $Q$ is not found to be in $\mathcal{L}$, then the volume of $Q$ is computed (if it has not already been computed), a new row is added to $\mathcal{L}$ and the volume and face vector are recorded.
 \begin{itemize}
   \item If the volume calculation fails for $Q$, this is also recorded.  If later a triangulation is found for which the volume calculation of $Q$ is successful, the row in $\mathcal{L}$ corresponding to $Q$ is updated with the information.
 \end{itemize}
\item[(6)] All polyhedra created which were found to be previously unrepresented in $\mathcal{L}$ are written from memory into an orb file.  The list $\mathcal{L}$ with its updated information is written to a file. 
\end{enumerate}

All of the programs mentioned above are freely available.  The {\tt python} progrms and the modified version of {\tt{Orb}}, along with the {\tt{Orb}} files I generated running the program, are freely available upon request.

\section{Composition}

I decided not to include a computer implementation of composition because of the computational cost of computing all possible compositions of two polyhedra.  To demonstrate this, note that every right-angled hyperbolic polyhedron has at least 12 pentagonal faces.  Unless some symmetry conditions can be exploited, the number of compositions between two polyhedra is at least $12 \times 12 \times 10 = 1440$ and this only counts compositions along pentagons.  There will be an even greater number if both polyhedra have hexagons, septagons, etc.  Finding the complete set of compositions of every pair of polyhedra in the existing list of smallest polyhedra would be very computationally expensive, and so my decision was to forgo it.


Many polyhedra which can be obtained by composition do not need composition to be built up from the L\"obell polyhedra -- edge-additions alone will suffice.  So such compositions will be included in the list as edge-children of the L\"obell polyhedra.  But this is not always the case.  An example of a polyhedron which cannot be obtained from edge-additions on the L\"obell polyhedra is the composition of the dodecahedron with itself denoted $L_5 \cup L_5$.  Note that by symmetry of $L_5$, there is only one possible composition of these two dodecahedra.  This composition can be thought of as doubling $L_5$ along any one of its faces.  Since this polyhedron fails to have any very good edges, decomposition is required to break it down into L\"obell polyhedra.  

I noted above that the program needs some initial set of polyhedra (seeds) to begin computing the list of smallest polyhedra.  Denote this initial set $\mathcal{I}$ and the union of $\mathcal{I}$ with set of all edge-descendants by $\mathcal{D}(\mathcal{I})$.  If $\mathcal{I}$ consists only of L\"obell polyhedra, then $\mathcal{D}(\mathcal{I})$ will not include $L_5 \cup L_5$.  Therefore, $L_5 \cup L_5$ is included in $\mathcal{I}$.  

So the question becomes, ``If $\mathcal{I}$ consists of the L\"obell polyhedra and $L_5 \cup L_5$, is $\mathcal{D}(\mathcal{I})$ the set of all right-angled hyperbolic polyhedra?"  The answer to this question is no.  This will be discussed more below.

I will prove the following weaker proposition:

\begin{theorem}\label{fifteen}
 If $\mathcal{I}$ consists of the L\"obell polyhedra and $L_5 \cup L_5$, then the set of right-angled hyperbolic polyhedra whose volume is less than 15 is contained in $\mathcal{D}(\mathcal{I})$.
\end{theorem}

\noindent \textbf{Proof:}  Corollary 2 says every right-angled polyhedron $P$ is obtainable by edge-additions and compositions of the L\"obell polyhedra and their edge-descendants.  By definition, if $P$ is obtainable from the L\"obell polyhedra using only edge-additions, then it is included in $\mathcal{D}(\mathcal{I})$.  In light of this, the only polyhedra which could be excluded from $\mathcal{D}(\mathcal{I})$ are compositions and their edge-descendants.  The theorem is proved if it is shown that all compositions which have volume less than 15 are contained in $\mathcal{D}(\mathcal{I})$.

Here is a simple result which is helpful to this end:

\begin{lemma} Suppose one of the following statements holds for a right-angled polyhedron $P$:
 \begin{enumerate}
  \item $P$ has a pair of very good edges $e_1$ and $e_2$ which satisfy the following conditions:  
  \begin{enumerate}
  \item $e_1$ and $e_2$ do not both belong to the same face, 
  \item if $e_1$ edge-connects faces $F_1$ to $F_2$ and $e_2$ edge-connects $G_1$ to $G_2$, then $F_i \neq G_j$ for $i,j = 1,2$,
  \item $e_1$ is not an edge of $G_1$ or $G_2$, and $e_2$ is not an edge of $F_1$ or $F_2$.
  \end{enumerate}
  \item $P$ has three very good edges, $e_{12}, e_{23}, e_{31}$, no two of which belong to the same face, and which edge-connect faces $F_1$ to $F_2$, $F_2$ to $F_3$, and $F_3$ to $F_1$.  
 \end{enumerate}
 Then any composition of polyhedra which involves $P$ contains a very good edge $e$ which has the property that edge-deleting $e$ will result in a polyhedron which is decomposable.
\end{lemma}

The utility of this result is that if $Q$ is the composition of two polyhedra $P_1$ and $P_2$, one of which satisfies the conditions of the lemma, then $Q$ must have a very good edge which when deleted results in a strictly smaller-volume composition.  If it has been shown that this smaller-volume composition is in our the list of edge-descendants $\mathcal{D}(\mathcal{I})$, then it follows that $Q$ will also be in $\mathcal{D}(\mathcal{I})$. 

\medskip

\noindent \textbf{Proof of Lemma:}  This proof will begin with a sub-lemma.

\begin{sub-lemma} Suppose $e$ is a very good edge of $P$.  Denote the faces of $P$ to which $e$ belongs $J_1$ and $J_2$, and the faces that $e$ edge-connects $K_1$ and $K_2$.  Denote by $eP$ the polyhedron which is obtained by edge-deletion of $e$.  Then the composition $P \cup Q$ along a face $F \subset P$ other than $K_1$, $K_2$, $J_1$, and $J_2$ has a very good edge which when edge-deleted, results in the composition $eP \cup Q$. 
\end{sub-lemma}

\noindent \textbf{Proof of Sub-Lemma:}  For ease of communication, it is useful to define topological inclusions of $P$ and $Q$ into the composition $P \cup Q$.  The polyhedron $P \cup Q$ has a distinguished prismatic $k$--circuit, call it $d$  corresponding to the boundary of the face $F$ that was the site of the composition.  The prismatic $k$--circuit bounds an embedded topological $k$--gon in $P \cup Q$ which is the site on the boundary of $P$ where $P$ and $Q$ were glued.  Splitting $P \cup Q$ along this topological $k$--gon produces two components which are topologically and combinatorially isomorphic to $P$ and $Q$.  This defines inclusions $i: P \to P\cup Q$ and $j: Q \to P\cup Q$, 

A simple observation about a composition involving $P$ is that it does not decrease the number of edges in any face of $P$ which is not the gluing site $F$. By this I mean that if $G$ is any face of $P$ other than $F$, $i(G)$ is a face of $P \cup Q$ which has at least as many edges as $G$.  Further, composition necessarily increases the number of edges of those faces which are adjacent to $F$.  By this I mean that if $C$ is a face of $P$ adjacent to $F$, then $i(C)$ is contained in a face of $P\cup Q$ which has more edges than $C$.

Therefore, returning to the context of the sub-lemma, since composition takes place at a face other than $K_1$ or $K_2$, these faces do not see a decrease in their number of edges.  Also, since the composition does not happen at $J_1$ or $J_2$, we know that $i(e)$ is an edge of $P \cup Q$.  Therefore, $i(e)$ edge-connects two large faces, namely $i(K_1)$ and $i(K_2)$.

For the purposes of establishing a contradiction, suppose that $i(e)$ is an edge intersected by a prismatic $5$--circuit $c$.  If $c$ does not intersect the distinguished prismatic $k$--circuit $d$ described above, then $c$ must be a prismatic 5--circuit contained entirely in $i(P)$.  This implies that $e$ is intersected by a prismatic 5--circuit in $P$. This is a contradiction.

So suppose $c$ intersects $d$.  Since $c$ and $d$ are simple closed loops on a sphere, they must intersect in at least two points.  These two points can get mapped back to $P$ by the inverse map $i^{-1}: i(P) \to P$.  Since $c$ is prismatic, these two points in $P$ lie on two different edges of the face $F$ which do not share an endpoint.  These two points can then be connected by a curve lying entirely within $F$ to connect the endpoints create a prismatic circuit in $P$.  Note that this prismatic circuit necessarily intersects fewer edges in $P$ than $c$ does in $P\cup Q$.  This is a contradiction.

This proves that $i(e)$ is a very good edge of $P \cup Q$.

Let $R$ denote the polyhedron which is obtained by edge-deleting $i(e)$.  Note that $d$ persists as a prismatic $k$--circuit in $R$ and splitting $R$ along the topological $k$--gon bounded by $d$ results in components combinatorially isomorphic to $eP$ and $Q$.  Therefore, $R$ is $eP \cup Q$.  

This ends the proof of the sub-lemma.
\qed

\medskip

Returning to the proof of the lemma, for any composition of $P$ and $Q$ along a face $F \subset P$, either of the conditions on $P$ described in the statement of the lemma guarantees the existence of a very good edge in $P$ which is disjoint from $F$.

Suppose $P$ satisfies the first condition of the lemma.  The faces which $e_1$ edge-connects are denoted $F_1$ and $F_2$ -- denote the faces to which $e_1$ belongs by $J_1$ and $J_2$.  If the face $F$ where the composition is to take place is not one of these four faces, then the conclusion of the lemma follows from the sub-lemma.

So suppose $F = F_1$.  Then by condition (b), $F$ is not $G_1$ nor $G_2$, and by condition (c), $e_2$ is not an edge of $F$.  Therefore, this choice of edge $e_2$ and gluing face $F$ satisfy the conditions of the sub-lemma, from which the lemma follows.  A similar argument works if $F = F_2$.

So suppose $F = J_1$.  Then $e_1$ is an edge of $F$.  Then by condition (a), $e_2$ is not an edge of $F$ and by condition (c), neither $G_1$ nor $G_2$ is $F$.  Thus, the lemma follows.

The above argument can be trivially adapted by relabeling indices to show that if $F = G_1$, $F = G_2$, or $e_2$ is an edge of $F$, then $e_1$ and $F$ satisfy the conditions of the sub-lemma.  Hence, the conclusion of the lemma follows from condition (1).

Now consider condition (2) of the lemma.  Suppose $F = F_1$.  Then $F$ is neither $F_2$ nor $F_3$.  Further, $e_{23}$ cannot be an edge of $F$ since trivalence would then imply $F_1$ is adjacent to both $F_2$ and $F_3$ which would preclude the existence of a very good edge (or any edge at all) between them.  Therefore, the choice of edge $e_{23}$ and choice of gluing site $F = F_1$ satisfy the conditions of the sub-lemma, and so the lemma follows.  This argument can be adapted to the cases where $F = F_2$ or $F = F_3$ as well.

Suppose $F$ has $e_{12}$ as an edge.  Then neither $e_{23}$ nor $e_{31}$ are edges of $F$.  Neither $F_1$ nor $F_2$ can be $F$ since $e_{12}$ edges connects $F_1$ and $F_2$.  Further, $F_3$ cannot be $F$ trivalence would imply that $F_3$ is adjacent to $F_1$ and $F_2$.  Therefore, the lemma follows.

This ends the proof of the lemma.
\qed

\medskip

Returning to the proof of the theorem, since only compositions of volume less than 15 need to be considered, I will restrict attention to compositions involving the 39 smallest-volume polyhedra in the set $\mathcal{D}(\mathcal{I})$ -- denote this subset $\mathcal{A}$. This is because any composition that involves a polyhedron in $\mathcal{D}(\mathcal{I})$ larger in volume than those in $\mathcal{A}$ is necessarily larger than 15 in volume \footnote{The 40th smallest polyhedron in $\mathcal{D}(\mathcal{I})$ has volume 10.72713 and so the smallest possible compoisition (which is with $L_5$) will have volume at least $10.72713 + 4.3062 > 15$}.   I will denote the $n$th smallest polyhedron in $\mathcal{A}$ by $A_n$.  So, for instance, $A_1 = L_5$, $A_2 = L_6$, and $A_3$ is the unique polyhedron obtained by edge-addition on $L_6$.

By manual observation of 1--skeleta, those polyhedra which satisfy the conditions of the lemma can be detected (see Appendix B for pictures of these 1--skeleta).  It happens that 13 of the polyhedra in $\mathcal{A}$ fail to satisfy the conditions of the lemma.  Denote this set of 13 polyhedra $\mathcal{B}$.  Interestingly,  $\mathcal{B}$ contains exactly those polyhedra which are not edge-descendants of $A_2 = L_6$.  Figure \ref{familytree} shows the family tree of edge-descendants.  Notice that $A_1$, $A_4$, $A_7$, $A_{11}$, $A_{12}$, $A_{18}$, $A_{23}$, $A_{26}$, $A_{29}$, $A_{30}$, $A_{35}$, $A_{39}$ are all missing from this family tree.  Adding $A_2$ completes the list $\mathcal{B}$.

Table 1 contains the volumes of all the elements of $\mathcal{A}$. Shaded in gray are those polyhedra which lie in $\mathcal{B}$.

\smallskip

\begin{center}
\begin{table}
\begin{tabular}{ | c | c | c |}
 \hline
 $\cellcolor{gray!25}\vol(A_1) = 4.30621$ &  $\vol(A_{14}) = 9.56488$ & $\vol(A_{27}) = 10.34848$\\
 \hline
 $\cellcolor{gray!25}\vol(A_2) = 6.02304$  &  $\vol(A_{15}) = 9.62756$ & $\vol(A_{28}) = 10.40429$\\ \hline 
 $\vol(A_3) = 6.96701$  &  $\vol(A_{16}) = 9.67726$ & $\cellcolor{gray!25}\vol(A_{29}) = 10.41604$\\ \hline 
 $\cellcolor{gray!25}\vol(A_4) = 7.56325$  &  $\vol(A_{17}) = 9.71117$ & $\cellcolor{gray!25}\vol(A_{30}) = 10.42605$\\ \hline 
 $\vol(A_5) = 7.86995$  &  $\cellcolor{gray!25}\vol(A_{18}) = 9.73084$ & $\vol(A_{31}) = 10.48885$\\ \hline 
 $\vol(A_6) = 8.00023$  &  $\vol(A_{19}) = 9.80384$ & $\vol(A_{32}) = 10.48909$\\ \hline 
 $\cellcolor{gray!25}\vol(A_7) = 8.61241$  &  $\vol(A_{20}) = 9.81423$ & $\vol(A_{33}) = 10.53439$\\ \hline 
 $\vol(A_8) = 8.67652$  &  $\vol(A_{21}) = 9.83568$ & $\vol(A_{34}) = 10.56155$\\ \hline 
 $\vol(A_9) = 8.86089$  &  $\vol(A_{22}) = 9.92355$ & \cellcolor{gray!25}$\vol(A_{35}) = 10.59201$\\ \hline 
 $\vol(A_{10}) = 8.9466$  &  $\cellcolor{gray!25}\vol(A_{23}) = 9.97683$ & $\vol(A_{36}) = 10.61998$\\ \hline 
 $\cellcolor{gray!25}\vol(A_{11}) =9.01905$  &  $\vol(A_{24}) = 9.97717$ & $\vol(A_{37}) = 10.63572$ \\ \hline 
 $\cellcolor{gray!25}\vol(A_{12}) =9.20916$  &  $\vol(A_{25}) = 10.21991$ & $\vol(A_{38}) = 10.67059$\\ \hline 
 $\vol(A_{13}) =9.47497$  &  $\cellcolor{gray!25}\vol(A_{26}) = 10.3378$ &\cellcolor{gray!25} $\vol(A_{39}) = 10.67059$\\ \hline 
\end{tabular}
\caption{This table contains the volumes of the 39 smallest right-angled hyperbolic polyhedra in $\mathcal{D}(\mathcal{I})$.  Shaded cells represent polyhedra which do not satisfy the conditions of the lemma.}
\end{table}
\end{center}

\smallskip

The task now is to show that compositions of two polyhedra in $\mathcal{B}$ and which have volume less than 15 are elements of $\mathcal{D}(\mathcal{I})$.  Once this has been shown, the theorem will be proved.

Here is the idea for how this proves the theorem.  Consider a composition of any two elements of $\mathcal{A}$, say $P \cup Q$.  If one of them satisfies the conditions of the lemma (that is, one of them is in $\mathcal{A} \setminus \mathcal{B}$), then an edge-deletion will obtain $eP \cup Q$.  Note that since the volume of $eP$ is strictly smaller than that of $P$, $eP$ is also in $\mathcal{A}$.  If either $eP$ or $Q$ satisfies the conditions of the lemma, then we can repeat this process.  Repeat until both of the components of the composition fail to satisfy the conditions of the lemma.  Then the composition is in $\mathcal{D}(\mathcal{I})$ as are all of its edge-descendants which includes $P \cup Q$.

Certain pairs of polyhedra do not need to considered since the volume of the composition is easily determined to be larger than 15 (recall the lower bound on composition $\vol(P \cup Q) \geq \vol(P) + \vol(Q)$).  For instance, for any polyhedron in $\mathcal{B}$ whose volume is larger than 9, only the composition with the dodecahedron $A_1$ needs to be considered since composing it with anything with larger volume will result in a polyhedron whose volume is larger than 15.  This reduces the amount of work to do considerably.

Pairs of polyhedra in $\mathcal{B}$ whose compositions need to be considered and shown to be in $\mathcal{D}(\mathcal{I})$ are: 
\begin{itemize}
\item $(A_1, *)$, where $*$ is any element of $\mathcal{B}$,
\item $(A_2, A_2), (A_2, A_4), (A_2, A_7)$
\end{itemize}
Note that the polyhedra $A_1$, $A_2$, $A_4$, $A_{11}$, and $A_{30}$ are the L\"obell polyhedra $L_5$, $L_6$, $L_7$, $L_8$ and $L_9$ respectively.  Also, $A_7$ is $L_5 \cup L_5$.

Many of these pairs have polyhedra with symmetry which reduces the number of compositions that need to be considered.  For example, the symmetries of the dodecahedron $A_1$ mean that the number of compositions between $A_1$ and any other polyhedron $P$ is determined by the number of pentagons that $P$ has as faces.  

To show that compositions of these elements of $\mathcal{B}$ are in $\mathcal{D}(\mathcal{I})$, I constructed the composition in {\tt{Orb}}'s graphical user interface (really, it is the double cover which I constructed).  This produces an orbifold which I can compare against the set $\mathcal{D}(\mathcal{I})$  (namely the polyhedra represented in the list $\mathcal{L}$) by checking for graph isomorphism.

\begin{figure}
\centering
\includegraphics[height=4cm, angle=20, origin=c]{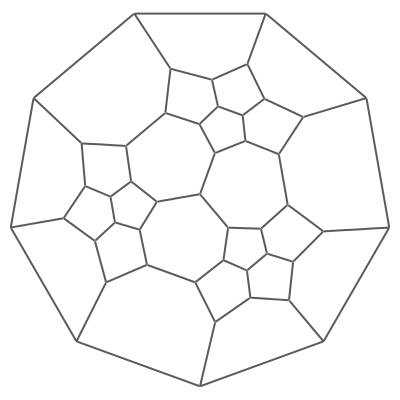}
\caption{This polyhedron has no very good edge.}
\label{novg}
\end{figure}

This method shows that every composition of pairs of elements of $\mathcal{B}$ is in $\mathcal{D}(\mathcal{I})$ except one!  This exceptional composition is $A_{39} \cup A_1$ along one of $A_{39}$'s fifteen pentagonal faces (see Figure \ref{novg}). Fortunately, this composition has volume 15.07032... and so does not contradict the theorem.  I will discuss this interesting counterexample more below.

This proves the theorem.
\qed
\smallskip

So to reiterate: every right-angled hyperbolic polyhedron with volume less than 15 can be obtained by repeated edge-additions on a L\"obell polyhedron or $L_5 \cup L_5$.  Therefore, the list of polyhedra whose volume is less than 15 can be obtained without having to implement composition algorithmically.  

\section{Results}

A table of volumes of the 825 smallest polyhedra appears in the appendix in this paper. The Tutte embeddings of the 100 smallest polyhedra also appears in the appendix.  In Figure \ref{scatterplot}, a scatter plot of the volumes of the first 200 smallest polyhedron is given.

\begin{figure}
\centering
 \includegraphics[trim={5cm 14cm 0 3cm}, clip]{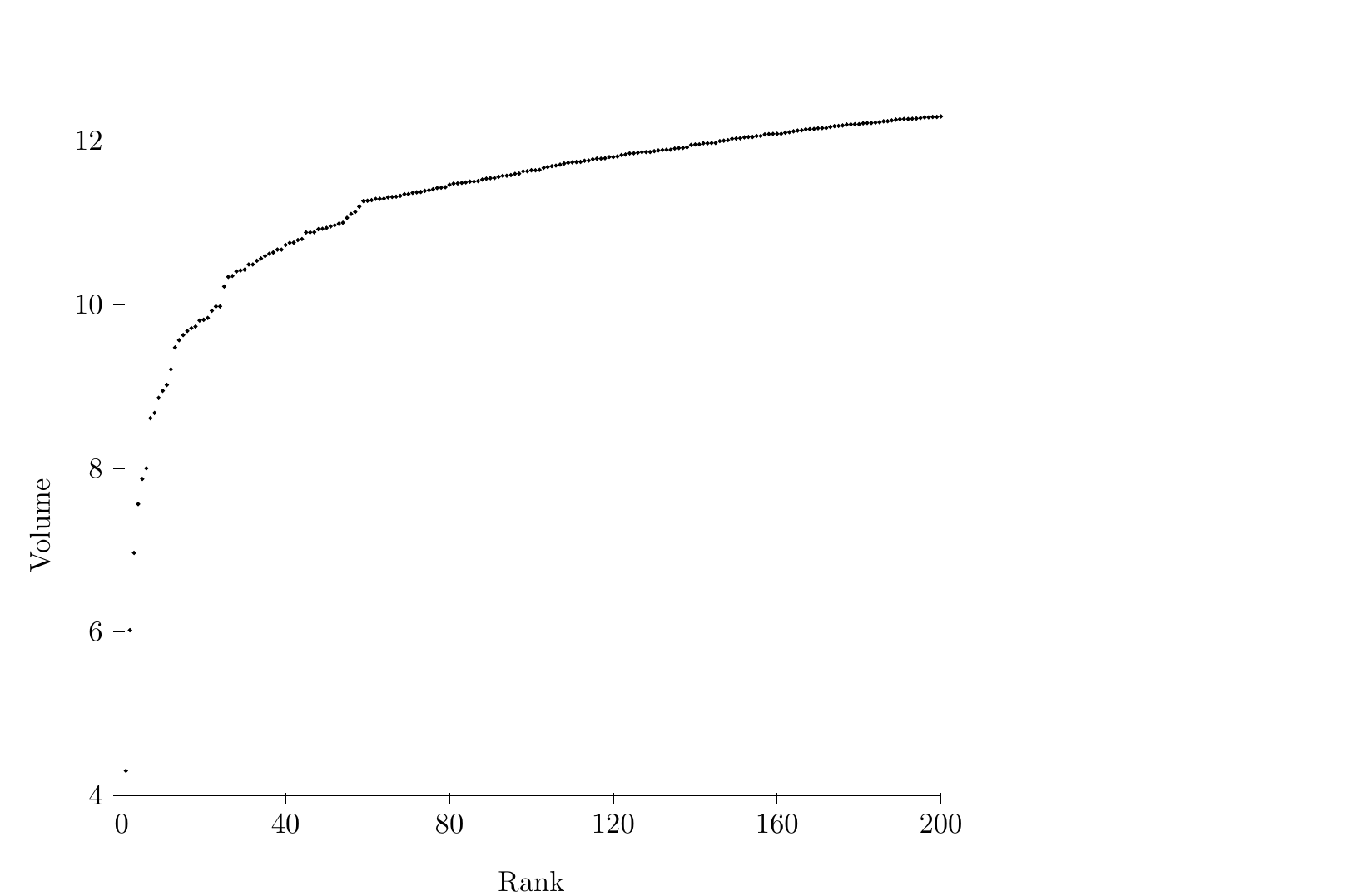}
 \caption{A scatter plot of the volumes of the first 200 smallest right-angled hyperbolic polyhedra.}
 \label{scatterplot}
\end{figure}

\section{Problems and Questions}
\begin{figure}
\centering
\includegraphics[height=4cm, angle = -77.14, origin=c]{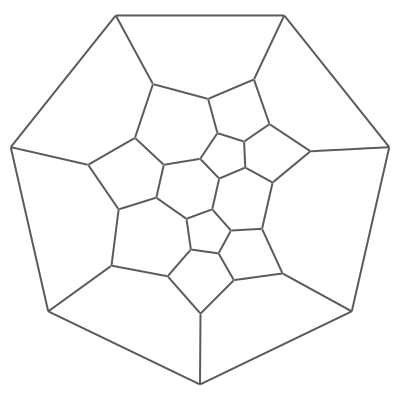}\includegraphics[height=4cm, angle=-67.5, origin=c]{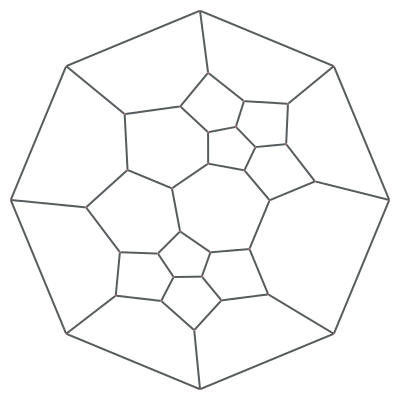}	
\caption{A pair of non-isometric right-angled polyhedra with the same volume (up to 12 digits). What I am calling $A_{38}$ is on the left, and $A_{39}$ is on the right.}
\label{nonisopair}
\end{figure}

Here are some questions which I think are interesting:

\begin{enumerate}
 \item The list as it currently stands consists of 825 polyhedra.  The 825th smallest right-angled polyhedron has volume 13.42033.  But Theorem \ref{fifteen} implies that this list can be greatly extended without needing to introduce composition. The obstacle is that {\tt Orb} is unable to compute the hyperbolic structure of certain right-angled hyperbolic polyhedra in my list and so I cannot accurately place them in the list.  If a way of obtaining a reliable volume for these polyhedra comes to exist, then the list could probably be extended to many thousands of polyhedra.
 
 \item The claim that every right-angled hyperbolic polyhedron is an edge-descendant of the L\"obell polyhedra or $L_5 \cup L_5$ is false -- a counterexample is provided by a composition of $A_{39}$ and $A_1 = L_5$ (see Figure \ref{novg}).  This polyhedron cannot be an edge-descendant of anything because it has no very good edges.  Of course, this polyhedron could be added to $\mathcal{I}$, and the list of smallest-volume polyhedra would then be extended without implementing composition, but there is no telling where the next such exceptional polyhedra will appear.  Is there a characterization of polyhedra which are not edge-descendants?
 
 \item $A_{39}$ is also notable because its {\tt Orb}-reported volume is equal to the {\tt Orb}-reported volume of $A_{38}$ up to twelve digits of precision.  This is not so unusual later in the list, when volumes become more dense.  But two polyhedra with almost the same volume occur so early in the list is remarkable and raises the question: Why?  Non-isometric polyhedra with the same volume can also arise by doubling a polyhedron across its various faces however $A_{38}$ and $A_{39}$ are not related in this way.  See Figure \ref{nonisopair}.  
 
 \item This list depends heavily on {\tt Orb} and {\tt graph-tool}. These are fantastic open-source programs but I lack the expertise to judge their numerical accuracy and correctness which are issues when working with very complicated polyhedra.  Are these programs reliable enough to trust all of the information found in this list?
 
 \item Being able to see these hyperbolic polyhedra in their native habitat (namely $\mathbb{H}^3$) would be beautiful. A nice project would be to write a program which can take the geometric triangulation data found in {\tt Orb} and create an image of the polyhedra in, say, the Poincar\'e ball model.  The program {\tt Geomview} is able to render such images, so this is really just a matter of translating from {\tt Orb}'s file format to {\tt Geomview}.
\end{enumerate}

\vspace{.3in}

\bibliographystyle{halpha}
\bibliography{ref}

\section*{Appendix A: Table of Volumes}

\noindent Rank:  Volume
\begin{multicols}{3}
\begin{tabbing}
Test \quad\= 99\quad\=\kill
1: \> 4.3062108 \\ 
2: \> 6.023046 \\ 
3: \> 6.9670114 \\ 
4: \> 7.5632491 \\ 
5: \> 7.8699479 \\ 
6: \> 8.0002343 \\ 
7: \> 8.6124152 \\ 
8: \> 8.6765244 \\ 
9: \> 8.8608973 \\ 
10: \> 8.9466057 \\ 
11: \> 9.0190528 \\ 
12: \> 9.2091595 \\ 
13: \> 9.4749695 \\ 
14: \> 9.5648792 \\ 
15: \> 9.6275647 \\ 
16: \> 9.6772650 \\ 
17: \> 9.7111700 \\ 
18: \> 9.7308471 \\ 
19: \> 9.8038472 \\ 
20: \> 9.8142328 \\ 
21: \> 9.8356794 \\ 
22: \> 9.9235512 \\ 
23: \> 9.9768361 \\ 
24: \> 9.9771698 \\ 
25: \> 10.2199075 \\ 
26: \> 10.3378015 \\ 
27: \> 10.3484815 \\ 
28: \> 10.4042904 \\ 
29: \> 10.4160438 \\ 
30: \> 10.4260522 \\ 
31: \> 10.4888519 \\ 
32: \> 10.489097 \\ 
33: \> 10.5343945 \\ 
34: \> 10.5615537 \\ 
35: \> 10.5920090 \\ 
36: \> 10.6199778 \\ 
37: \> 10.6357266 \\ 
38: \> 10.6705888 \\ 
39: \> 10.6705888 \\ 
40: \> 10.7271348 \\ 
41: \> 10.7534763 \\ 
42: \> 10.7563401 \\ 
43: \> 10.7874579 \\ 
44: \> 10.7990421 \\ 
45: \> 10.8804124 \\ 
46: \> 10.8816862 \\ 
47: \> 10.8835451 \\ 
48: \> 10.9208994 \\ 
49: \> 10.9248751 \\ 
50: \> 10.9363484 \\ 
51: \> 10.9543051 \\ 
52: \> 10.9682397 \\ 
53: \> 10.9860570 \\ 
54: \> 10.9996896 \\ 
55: \> 11.0587626 \\ 
56: \> 11.1061236 \\ 
57: \> 11.1302424 \\ 
58: \> 11.1955201 \\ 
59: \> 11.263206 \\ 
60: \> 11.2679748 \\ 
61: \> 11.2751917 \\ 
62: \> 11.2904576 \\ 
63: \> 11.2913353 \\ 
64: \> 11.2932291 \\ 
65: \> 11.3094811 \\ 
66: \> 11.3142710 \\ 
67: \> 11.3182574 \\ 
68: \> 11.3281616 \\ 
69: \> 11.3480005 \\ 
70: \> 11.3502332 \\ 
71: \> 11.3639865 \\ 
72: \> 11.3709239 \\ 
73: \> 11.3743487 \\ 
74: \> 11.3888812 \\ 
75: \> 11.3959620 \\ 
76: \> 11.4070515 \\ 
77: \> 11.4229627 \\ 
78: \> 11.4275803 \\ 
79: \> 11.4335963 \\ 
80: \> 11.4635122 \\ 
81: \> 11.4778773 \\ 
82: \> 11.4795650 \\ 
83: \> 11.4865115 \\ 
84: \> 11.4917848 \\ 
85: \> 11.5019129 \\ 
86: \> 11.5026744 \\ 
87: \> 11.5079212 \\ 
88: \> 11.5262055 \\ 
89: \> 11.5375492 \\ 
90: \> 11.5430117 \\ 
91: \> 11.545118 \\ 
92: \> 11.5594844 \\ 
93: \> 11.5719948 \\ 
94: \> 11.5732855 \\ 
95: \> 11.581167 \\ 
96: \> 11.5962092 \\ 
97: \> 11.5994297 \\ 
98: \> 11.6272353 \\ 
99: \> 11.6282329 \\ 
100: \> 11.6394320 \\ 
101: \> 11.639562 \\ 
102: \> 11.6445996 \\ 
103: \> 11.6696888 \\ 
104: \> 11.6794851 \\ 
105: \> 11.6904395 \\ 
106: \> 11.6970789 \\ 
107: \> 11.7084617 \\ 
108: \> 11.7220185 \\ 
109: \> 11.7318198 \\ 
110: \> 11.7371977 \\ 
111: \> 11.7408513 \\ 
112: \> 11.7423134 \\ 
113: \> 11.7559468 \\ 
114: \> 11.7578225 \\ 
115: \> 11.7751062 \\ 
116: \> 11.7816162 \\ 
117: \> 11.7823418 \\ 
118: \> 11.787139 \\ 
119: \> 11.80131 \\ 
120: \> 11.8013150 \\ 
121: \> 11.8091731 \\ 
122: \> 11.8262458 \\ 
123: \> 11.8316258 \\ 
124: \> 11.8463446 \\ 
125: \> 11.8477745 \\ 
126: \> 11.8539277 \\ 
127: \> 11.8612950 \\ 
128: \> 11.861906 \\ 
129: \> 11.8626576 \\ 
130: \> 11.873384 \\ 
131: \> 11.8816873 \\ 
132: \> 11.8878994 \\ 
133: \> 11.8911801 \\ 
134: \> 11.8911801 \\ 
135: \> 11.9051986 \\ 
136: \> 11.9111339 \\ 
137: \> 11.9130597 \\ 
138: \> 11.9201173 \\ 
139: \> 11.9496613 \\ 
140: \> 11.9548014 \\ 
141: \> 11.9562534 \\ 
142: \> 11.9693268 \\ 
143: \> 11.9695181 \\ 
144: \> 11.9722376 \\ 
145: \> 11.9734938 \\ 
146: \> 11.9958248 \\ 
147: \> 12.0009332 \\ 
148: \> 12.0071733 \\ 
149: \> 12.0251598 \\ 
150: \> 12.0284546 \\ 
151: \> 12.0307865 \\ 
152: \> 12.0421556 \\ 
153: \> 12.0460920 \\ 
154: \> 12.0460920 \\ 
155: \> 12.0577601 \\ 
156: \> 12.0588613 \\ 
157: \> 12.0774156 \\ 
158: \> 12.0816685 \\ 
159: \> 12.0841913 \\ 
160: \> 12.0841913 \\ 
161: \> 12.0855209 \\ 
162: \> 12.0984817 \\ 
163: \> 12.1038544 \\ 
164: \> 12.1147206 \\ 
165: \> 12.1233416 \\ 
166: \> 12.1269916 \\ 
167: \> 12.1405768 \\ 
168: \> 12.1415323 \\ 
169: \> 12.1444998 \\ 
170: \> 12.1517046 \\ 
171: \> 12.1545352 \\ 
172: \> 12.1553815 \\ 
173: \> 12.1682871 \\ 
174: \> 12.1781926 \\ 
175: \> 12.1810448 \\ 
176: \> 12.1863417 \\ 
177: \> 12.1976346 \\ 
178: \> 12.2001378 \\ 
179: \> 12.2011148 \\ 
180: \> 12.2017855 \\ 
181: \> 12.2124156 \\ 
182: \> 12.2156059 \\ 
183: \> 12.2160447 \\ 
184: \> 12.2211128 \\ 
185: \> 12.2249334 \\ 
186: \> 12.2372753 \\ 
187: \> 12.2380728 \\ 
188: \> 12.2478996 \\ 
189: \> 12.2568250 \\ 
190: \> 12.2638732 \\ 
191: \> 12.2654041 \\ 
192: \> 12.2654041 \\ 
193: \> 12.2677948 \\ 
194: \> 12.2717135 \\ 
195: \> 12.276929 \\ 
196: \> 12.2849205 \\ 
197: \> 12.2854781 \\ 
198: \> 12.2900996 \\ 
199: \> 12.2907040 \\ 
200: \> 12.2960383 \\ 
201: \> 12.3030129 \\ 
202: \> 12.3085011 \\ 
203: \> 12.3190359 \\ 
204: \> 12.329375 \\ 
205: \> 12.3320289 \\ 
206: \> 12.3384939 \\ 
207: \> 12.3401788 \\ 
208: \> 12.3457956 \\ 
209: \> 12.3522024 \\ 
210: \> 12.3528353 \\ 
211: \> 12.3536105 \\ 
212: \> 12.3554360 \\ 
213: \> 12.3642145 \\ 
214: \> 12.3642358 \\ 
215: \> 12.3711438 \\ 
216: \> 12.3715722 \\ 
217: \> 12.3803662 \\ 
218: \> 12.3806750 \\ 
219: \> 12.3817306 \\ 
220: \> 12.3817888 \\ 
221: \> 12.3839007 \\ 
222: \> 12.3842394 \\ 
223: \> 12.3887363 \\ 
224: \> 12.3974717 \\ 
225: \> 12.3990550 \\ 
226: \> 12.4033673 \\ 
227: \> 12.4048315 \\ 
228: \> 12.4115128 \\ 
229: \> 12.4130114 \\ 
230: \> 12.4159192 \\ 
231: \> 12.4168618 \\ 
232: \> 12.4272231 \\ 
233: \> 12.4277649 \\ 
234: \> 12.4294935 \\ 
235: \> 12.430078 \\ 
236: \> 12.4309411 \\ 
237: \> 12.4342336 \\ 
238: \> 12.4377764 \\ 
239: \> 12.4402367 \\ 
240: \> 12.4410680 \\ 
241: \> 12.4447933 \\ 
242: \> 12.4468324 \\ 
243: \> 12.451627 \\ 
244: \> 12.4535759 \\ 
245: \> 12.4555261 \\ 
246: \> 12.4670205 \\ 
247: \> 12.4679602 \\ 
248: \> 12.4807071 \\ 
249: \> 12.4825054 \\ 
250: \> 12.4859467 \\ 
251: \> 12.4903569 \\ 
252: \> 12.4920523 \\ 
253: \> 12.4926391 \\ 
254: \> 12.4963117 \\ 
255: \> 12.503454 \\ 
256: \> 12.5048058 \\ 
257: \> 12.5059314 \\ 
258: \> 12.5080834 \\ 
259: \> 12.5131029 \\ 
260: \> 12.513552 \\ 
261: \> 12.5146366 \\ 
262: \> 12.5157225 \\ 
263: \> 12.5211267 \\ 
264: \> 12.5230388 \\ 
265: \> 12.526865 \\ 
266: \> 12.5310112 \\ 
267: \> 12.5343536 \\ 
268: \> 12.5346198 \\ 
269: \> 12.5355372 \\ 
270: \> 12.5362290 \\ 
271: \> 12.5370769 \\ 
272: \> 12.5383370 \\ 
273: \> 12.5406134 \\ 
274: \> 12.5421341 \\ 
275: \> 12.5480161 \\ 
276: \> 12.5497054 \\ 
277: \> 12.5559033 \\ 
278: \> 12.5600217 \\ 
279: \> 12.5693664 \\ 
280: \> 12.5709272 \\ 
281: \> 12.5719055 \\ 
282: \> 12.5770091 \\ 
283: \> 12.5797398 \\ 
284: \> 12.5883592 \\ 
285: \> 12.5916877 \\ 
286: \> 12.5919487 \\ 
287: \> 12.5988561 \\ 
288: \> 12.6001356 \\ 
289: \> 12.6119078 \\ 
290: \> 12.6127051 \\ 
291: \> 12.6133512 \\ 
292: \> 12.6146167 \\ 
293: \> 12.6147991 \\ 
294: \> 12.6159185 \\ 
295: \> 12.6174108 \\ 
296: \> 12.6180672 \\ 
297: \> 12.6207832 \\ 
298: \> 12.6246110 \\ 
299: \> 12.6247937 \\ 
300: \> 12.6270088 \\ 
301: \> 12.630289 \\ 
302: \> 12.6346334 \\ 
303: \> 12.6372939 \\ 
304: \> 12.646941 \\ 
305: \> 12.6479120 \\ 
306: \> 12.6500163 \\ 
307: \> 12.6514822 \\ 
308: \> 12.6587897 \\ 
309: \> 12.6606309 \\ 
310: \> 12.6613411 \\ 
311: \> 12.6661597 \\ 
312: \> 12.6664503 \\ 
313: \> 12.6670277 \\ 
314: \> 12.6672765 \\ 
315: \> 12.6687087 \\ 
316: \> 12.6693988 \\ 
317: \> 12.6699947 \\ 
318: \> 12.6713641 \\ 
319: \> 12.6742252 \\ 
320: \> 12.6753857 \\ 
321: \> 12.6794010 \\ 
322: \> 12.6806399 \\ 
323: \> 12.6828872 \\ 
324: \> 12.6832494 \\ 
325: \> 12.6860769 \\ 
326: \> 12.6871951 \\ 
327: \> 12.697821 \\ 
328: \> 12.6985894 \\ 
329: \> 12.7025794 \\ 
330: \> 12.7033136 \\ 
331: \> 12.7072758 \\ 
332: \> 12.7087943 \\ 
333: \> 12.7089644 \\ 
334: \> 12.7097704 \\ 
335: \> 12.7112824 \\ 
336: \> 12.7134298 \\ 
337: \> 12.7157367 \\ 
338: \> 12.7195313 \\ 
339: \> 12.7214472 \\ 
340: \> 12.7232730 \\ 
341: \> 12.7262254 \\ 
342: \> 12.7294877 \\ 
343: \> 12.7330111 \\ 
344: \> 12.7347908 \\ 
345: \> 12.7384981 \\ 
346: \> 12.7399561 \\ 
347: \> 12.7401153 \\ 
348: \> 12.7406675 \\ 
349: \> 12.7432168 \\ 
350: \> 12.7436003 \\ 
351: \> 12.7442571 \\ 
352: \> 12.7455923 \\ 
353: \> 12.7455926 \\ 
354: \> 12.7460513 \\ 
355: \> 12.7460515 \\ 
356: \> 12.7482225 \\ 
357: \> 12.7585029 \\ 
358: \> 12.7599093 \\ 
359: \> 12.7645973 \\ 
360: \> 12.7661321 \\ 
361: \> 12.7667392 \\ 
362: \> 12.7672454 \\ 
363: \> 12.7677000 \\ 
364: \> 12.7685012 \\ 
365: \> 12.7693876 \\ 
366: \> 12.776313 \\ 
367: \> 12.7769034 \\ 
368: \> 12.7820728 \\ 
369: \> 12.7822353 \\ 
370: \> 12.7831954 \\ 
371: \> 12.7833781 \\ 
372: \> 12.7857129 \\ 
373: \> 12.7857129 \\ 
374: \> 12.7920778 \\ 
375: \> 12.7925214 \\ 
376: \> 12.7948392 \\ 
377: \> 12.7979026 \\ 
378: \> 12.7996480 \\ 
379: \> 12.8016266 \\ 
380: \> 12.8018956 \\ 
381: \> 12.8045382 \\ 
382: \> 12.8066760 \\ 
383: \> 12.8107155 \\ 
384: \> 12.8174862 \\ 
385: \> 12.8183900 \\ 
386: \> 12.821584 \\ 
387: \> 12.8221324 \\ 
388: \> 12.8226626 \\ 
389: \> 12.8231303 \\ 
390: \> 12.8245185 \\ 
391: \> 12.8256219 \\ 
392: \> 12.8272844 \\ 
393: \> 12.8282458 \\ 
394: \> 12.8284383 \\ 
395: \> 12.8363006 \\ 
396: \> 12.8365463 \\ 
397: \> 12.8416379 \\ 
398: \> 12.8433256 \\ 
399: \> 12.8441243 \\ 
400: \> 12.8447142 \\ 
401: \> 12.8476691 \\ 
402: \> 12.8493546 \\ 
403: \> 12.8607115 \\ 
404: \> 12.8609814 \\ 
405: \> 12.8661875 \\ 
406: \> 12.8708785 \\ 
407: \> 12.8718074 \\ 
408: \> 12.8718715 \\ 
409: \> 12.8726247 \\ 
410: \> 12.8749811 \\ 
411: \> 12.8756532 \\ 
412: \> 12.8763805 \\ 
413: \> 12.8765835 \\ 
414: \> 12.8766261 \\ 
415: \> 12.8789947 \\ 
416: \> 12.8792251 \\ 
417: \> 12.8808630 \\ 
418: \> 12.8815061 \\ 
419: \> 12.8843508 \\ 
420: \> 12.8885622 \\ 
421: \> 12.8927067 \\ 
422: \> 12.893312 \\ 
423: \> 12.8934687 \\ 
424: \> 12.8960354 \\ 
425: \> 12.9022597 \\ 
426: \> 12.9056424 \\ 
427: \> 12.9077398 \\ 
428: \> 12.9101791 \\ 
429: \> 12.9114325 \\ 
430: \> 12.9124632 \\ 
431: \> 12.9127128 \\ 
432: \> 12.9139315 \\ 
433: \> 12.914183 \\ 
434: \> 12.9161902 \\ 
435: \> 12.916207 \\ 
436: \> 12.9186228 \\ 
437: \> 12.9186228 \\ 
438: \> 12.9205532 \\ 
439: \> 12.9246735 \\ 
440: \> 12.9250878 \\ 
441: \> 12.9254826 \\ 
442: \> 12.9316159 \\ 
443: \> 12.9320131 \\ 
444: \> 12.9332599 \\ 
445: \> 12.9342855 \\ 
446: \> 12.9379678 \\ 
447: \> 12.9380414 \\ 
448: \> 12.9482748 \\ 
449: \> 12.9507107 \\ 
450: \> 12.952697 \\ 
451: \> 12.9542260 \\ 
452: \> 12.9564114 \\ 
453: \> 12.9574588 \\ 
454: \> 12.9594551 \\ 
455: \> 12.9621565 \\ 
456: \> 12.9623225 \\ 
457: \> 12.9675033 \\ 
458: \> 12.9727244 \\ 
459: \> 12.9769842 \\ 
460: \> 12.9787127 \\ 
461: \> 12.9787488 \\ 
462: \> 12.978926 \\ 
463: \> 12.9803944 \\ 
464: \> 12.9810244 \\ 
465: \> 12.9821528 \\ 
466: \> 12.9826962 \\ 
467: \> 12.9841541 \\ 
468: \> 12.9852469 \\ 
469: \> 12.9870666 \\ 
470: \> 12.9871326 \\ 
471: \> 12.9875988 \\ 
472: \> 12.9892227 \\ 
473: \> 12.9894036 \\ 
474: \> 12.9898039 \\ 
475: \> 12.9924142 \\ 
476: \> 12.9936663 \\ 
477: \> 12.9941005 \\ 
478: \> 12.9943858 \\ 
479: \> 12.9951085 \\ 
480: \> 12.9952274 \\ 
481: \> 12.9961182 \\ 
482: \> 12.9974405 \\ 
483: \> 12.9979937 \\ 
484: \> 13.0009026 \\ 
485: \> 13.0011358 \\ 
486: \> 13.0012106 \\ 
487: \> 13.0040128 \\ 
488: \> 13.0054875 \\ 
489: \> 13.0072294 \\ 
490: \> 13.0074697 \\ 
491: \> 13.0084985 \\ 
492: \> 13.0090672 \\ 
493: \> 13.0114696 \\ 
494: \> 13.0128611 \\ 
495: \> 13.0147794 \\ 
496: \> 13.0151189 \\ 
497: \> 13.0180899 \\ 
498: \> 13.0181706 \\ 
499: \> 13.0202361 \\ 
500: \> 13.020657 \\ 
501: \> 13.0224109 \\ 
502: \> 13.0234614 \\ 
503: \> 13.0266063 \\ 
504: \> 13.0269017 \\ 
505: \> 13.0295852 \\ 
506: \> 13.0296075 \\ 
507: \> 13.032234 \\ 
508: \> 13.0380346 \\ 
509: \> 13.0409708 \\ 
510: \> 13.0418263 \\ 
511: \> 13.0456395 \\ 
512: \> 13.0459132 \\ 
513: \> 13.0543813 \\ 
514: \> 13.0572175 \\ 
515: \> 13.0573782 \\ 
516: \> 13.0617630 \\ 
517: \> 13.0619967 \\ 
518: \> 13.0636692 \\ 
519: \> 13.063861 \\ 
520: \> 13.0669990 \\ 
521: \> 13.0670094 \\ 
522: \> 13.0696241 \\ 
523: \> 13.0705683 \\ 
524: \> 13.0706272 \\ 
525: \> 13.0708697 \\ 
526: \> 13.0711347 \\ 
527: \> 13.0712706 \\ 
528: \> 13.0779935 \\ 
529: \> 13.0783591 \\ 
530: \> 13.0801169 \\ 
531: \> 13.0862297 \\ 
532: \> 13.0874702 \\ 
533: \> 13.0896096 \\ 
534: \> 13.0923100 \\ 
535: \> 13.0927300 \\ 
536: \> 13.0931814 \\ 
537: \> 13.0934364 \\ 
538: \> 13.0936906 \\ 
539: \> 13.0967306 \\ 
540: \> 13.0977768 \\ 
541: \> 13.1030816 \\ 
542: \> 13.1040997 \\ 
543: \> 13.1046457 \\ 
544: \> 13.1049271 \\ 
545: \> 13.1053328 \\ 
546: \> 13.1062826 \\ 
547: \> 13.1097041 \\ 
548: \> 13.1111005 \\ 
549: \> 13.1130951 \\ 
550: \> 13.1140612 \\ 
551: \> 13.1162075 \\ 
552: \> 13.1169051 \\ 
553: \> 13.1174496 \\ 
554: \> 13.1186799 \\ 
555: \> 13.1210031 \\ 
556: \> 13.1213456 \\ 
557: \> 13.1222987 \\ 
558: \> 13.1231679 \\ 
559: \> 13.1247713 \\ 
560: \> 13.1268511 \\ 
561: \> 13.1292634 \\ 
562: \> 13.1304211 \\ 
563: \> 13.1314807 \\ 
564: \> 13.1317051 \\ 
565: \> 13.1343683 \\ 
566: \> 13.1357657 \\ 
567: \> 13.1382809 \\ 
568: \> 13.1388930 \\ 
569: \> 13.1390365 \\ 
570: \> 13.1444814 \\ 
571: \> 13.1456676 \\ 
572: \> 13.1460156 \\ 
573: \> 13.1470214 \\ 
574: \> 13.1477861 \\ 
575: \> 13.1489491 \\ 
576: \> 13.1522210 \\ 
577: \> 13.1555327 \\ 
578: \> 13.15596 \\ 
579: \> 13.1600088 \\ 
580: \> 13.1600749 \\ 
581: \> 13.1619423 \\ 
582: \> 13.1621361 \\ 
583: \> 13.1640136 \\ 
584: \> 13.1650307 \\ 
585: \> 13.1652047 \\ 
586: \> 13.1665945 \\ 
587: \> 13.1679649 \\ 
588: \> 13.168069 \\ 
589: \> 13.1690101 \\ 
590: \> 13.1702350 \\ 
591: \> 13.1715685 \\ 
592: \> 13.1721263 \\ 
593: \> 13.1738399 \\ 
594: \> 13.1739971 \\ 
595: \> 13.1756584 \\ 
596: \> 13.1757761 \\ 
597: \> 13.1760025 \\ 
598: \> 13.1786547 \\ 
599: \> 13.1817410 \\ 
600: \> 13.1832403 \\ 
601: \> 13.183777 \\ 
602: \> 13.1856221 \\ 
603: \> 13.1862884 \\ 
604: \> 13.1875498 \\ 
605: \> 13.1897624 \\ 
606: \> 13.1897771 \\ 
607: \> 13.1909487 \\ 
608: \> 13.1926681 \\ 
609: \> 13.1929176 \\ 
610: \> 13.1936068 \\ 
611: \> 13.1942405 \\ 
612: \> 13.1945769 \\ 
613: \> 13.1965696 \\ 
614: \> 13.2021880 \\ 
615: \> 13.2027053 \\ 
616: \> 13.2027089 \\ 
617: \> 13.2037099 \\ 
618: \> 13.2060858 \\ 
619: \> 13.2065159 \\ 
620: \> 13.2073974 \\ 
621: \> 13.2083242 \\ 
622: \> 13.2104582 \\ 
623: \> 13.2129831 \\ 
624: \> 13.2170375 \\ 
625: \> 13.2189832 \\ 
626: \> 13.2191858 \\ 
627: \> 13.2192984 \\ 
628: \> 13.2210006 \\ 
629: \> 13.2214908 \\ 
630: \> 13.2215990 \\ 
631: \> 13.2226254 \\ 
632: \> 13.2237937 \\ 
633: \> 13.2242809 \\ 
634: \> 13.2254565 \\ 
635: \> 13.2310355 \\ 
636: \> 13.23436 \\ 
637: \> 13.2363318 \\ 
638: \> 13.2382701 \\ 
639: \> 13.2387091 \\ 
640: \> 13.238879 \\ 
641: \> 13.2409922 \\ 
642: \> 13.241284 \\ 
643: \> 13.2424799 \\ 
644: \> 13.2427474 \\ 
645: \> 13.2430734 \\ 
646: \> 13.2437962 \\ 
647: \> 13.2443978 \\ 
648: \> 13.2444270 \\ 
649: \> 13.2458608 \\ 
650: \> 13.2461749 \\ 
651: \> 13.2461749 \\ 
652: \> 13.2491505 \\ 
653: \> 13.2494142 \\ 
654: \> 13.2507197 \\ 
655: \> 13.2543603 \\ 
656: \> 13.2546591 \\ 
657: \> 13.2547085 \\ 
658: \> 13.2551004 \\ 
659: \> 13.2555559 \\ 
660: \> 13.2558775 \\ 
661: \> 13.2566728 \\ 
662: \> 13.2580638 \\ 
663: \> 13.2590908 \\ 
664: \> 13.262002 \\ 
665: \> 13.2620194 \\ 
666: \> 13.2623329 \\ 
667: \> 13.2627342 \\ 
668: \> 13.264521 \\ 
669: \> 13.2655601 \\ 
670: \> 13.268484 \\ 
671: \> 13.2688609 \\ 
672: \> 13.2720985 \\ 
673: \> 13.2723909 \\ 
674: \> 13.2723909 \\ 
675: \> 13.2724289 \\ 
676: \> 13.2731006 \\ 
677: \> 13.2736291 \\ 
678: \> 13.2737198 \\ 
679: \> 13.2737812 \\ 
680: \> 13.2746987 \\ 
681: \> 13.2755816 \\ 
682: \> 13.2775572 \\ 
683: \> 13.2784161 \\ 
684: \> 13.2787393 \\ 
685: \> 13.2793568 \\ 
686: \> 13.2802927 \\ 
687: \> 13.2827828 \\ 
688: \> 13.2834887 \\ 
689: \> 13.2871770 \\ 
690: \> 13.2876890 \\ 
691: \> 13.2884589 \\ 
692: \> 13.2896787 \\ 
693: \> 13.2921035 \\ 
694: \> 13.2937144 \\ 
695: \> 13.3003419 \\ 
696: \> 13.3009378 \\ 
697: \> 13.3031293 \\ 
698: \> 13.3034265 \\ 
699: \> 13.3047731 \\ 
700: \> 13.3052398 \\ 
701: \> 13.3068856 \\ 
702: \> 13.3079604 \\ 
703: \> 13.3091544 \\ 
704: \> 13.31075 \\ 
705: \> 13.3107847 \\ 
706: \> 13.3109185 \\ 
707: \> 13.3136447 \\ 
708: \> 13.3138381 \\ 
709: \> 13.3163635 \\ 
710: \> 13.3179858 \\ 
711: \> 13.3183921 \\ 
712: \> 13.3185439 \\ 
713: \> 13.318907 \\ 
714: \> 13.3191095 \\ 
715: \> 13.3191284 \\ 
716: \> 13.319623 \\ 
717: \> 13.3202038 \\ 
718: \> 13.3212566 \\ 
719: \> 13.3223716 \\ 
720: \> 13.3239867 \\ 
721: \> 13.3259183 \\ 
722: \> 13.3266908 \\ 
723: \> 13.3278201 \\ 
724: \> 13.3279098 \\ 
725: \> 13.3283104 \\ 
726: \> 13.3284008 \\ 
727: \> 13.3296515 \\ 
728: \> 13.3315609 \\ 
729: \> 13.3321591 \\ 
730: \> 13.3327367 \\ 
731: \> 13.3341648 \\ 
732: \> 13.3355411 \\ 
733: \> 13.3388411 \\ 
734: \> 13.339317 \\ 
735: \> 13.3393506 \\ 
736: \> 13.3400671 \\ 
737: \> 13.3407227 \\ 
738: \> 13.3408317 \\ 
739: \> 13.3422229 \\ 
740: \> 13.3428668 \\ 
741: \> 13.3453554 \\ 
742: \> 13.3454021 \\ 
743: \> 13.3470733 \\ 
744: \> 13.3474075 \\ 
745: \> 13.3475313 \\ 
746: \> 13.3480883 \\ 
747: \> 13.3493078 \\ 
748: \> 13.3497468 \\ 
749: \> 13.3502204 \\ 
750: \> 13.3509879 \\ 
751: \> 13.3513294 \\ 
752: \> 13.3516647 \\ 
753: \> 13.3524406 \\ 
754: \> 13.3536589 \\ 
755: \> 13.3558375 \\ 
756: \> 13.3571361 \\ 
757: \> 13.357809 \\ 
758: \> 13.3584357 \\ 
759: \> 13.3586095 \\ 
760: \> 13.3591299 \\ 
761: \> 13.3599534 \\ 
762: \> 13.3599811 \\ 
763: \> 13.3616158 \\ 
764: \> 13.3618647 \\ 
765: \> 13.362095 \\ 
766: \> 13.3637717 \\ 
767: \> 13.3663303 \\ 
768: \> 13.3667382 \\ 
769: \> 13.3677394 \\ 
770: \> 13.3682944 \\ 
771: \> 13.3704383 \\ 
772: \> 13.3707603 \\ 
773: \> 13.3707862 \\ 
774: \> 13.373108 \\ 
775: \> 13.374108 \\ 
776: \> 13.3745796 \\ 
777: \> 13.3767745 \\ 
778: \> 13.3783487 \\ 
779: \> 13.3783487 \\ 
780: \> 13.3836986 \\ 
781: \> 13.3853641 \\ 
782: \> 13.3862236 \\ 
783: \> 13.3864672 \\ 
784: \> 13.3875324 \\ 
785: \> 13.3877963 \\ 
786: \> 13.3883749 \\ 
787: \> 13.3883982 \\ 
788: \> 13.3892607 \\ 
789: \> 13.3903343 \\ 
790: \> 13.3903839 \\ 
791: \> 13.3914604 \\ 
792: \> 13.3921628 \\ 
793: \> 13.3925619 \\ 
794: \> 13.393811 \\ 
795: \> 13.3958043 \\ 
796: \> 13.3967555 \\ 
797: \> 13.3974018 \\ 
798: \> 13.3984599 \\ 
799: \> 13.3996204 \\ 
800: \> 13.4001767 \\ 
801: \> 13.4004630 \\ 
802: \> 13.4009408 \\ 
803: \> 13.4041059 \\ 
804: \> 13.4048751 \\ 
805: \> 13.4053031 \\ 
806: \> 13.4056015 \\ 
807: \> 13.4057424 \\ 
808: \> 13.4086268 \\ 
809: \> 13.4088302 \\ 
810: \> 13.4099862 \\ 
811: \> 13.4112902 \\ 
812: \> 13.4116875 \\ 
813: \> 13.4135964 \\ 
814: \> 13.4138043 \\ 
815: \> 13.4142589 \\ 
816: \> 13.4155186 \\ 
817: \> 13.4156955 \\ 
818: \> 13.4157039 \\ 
819: \> 13.416979 \\ 
820: \> 13.4172534 \\ 
821: \> 13.4178439 \\ 
822: \> 13.4181377 \\ 
823: \> 13.4197214 \\ 
824: \> 13.420017 \\ 
825: \> 13.42033 \\ 
\end{tabbing}
\end{multicols}

\vspace{.1in}

\section*{Appendix B: 1-skeleta of the 100 smallest polyhedra}
\begin{tabular}{|m{0.22\textwidth} | m{0.22\textwidth} | m{0.22\textwidth} |  m{0.22\textwidth} | } 
\hline 
1
\includegraphics[height=3cm, origin=c]{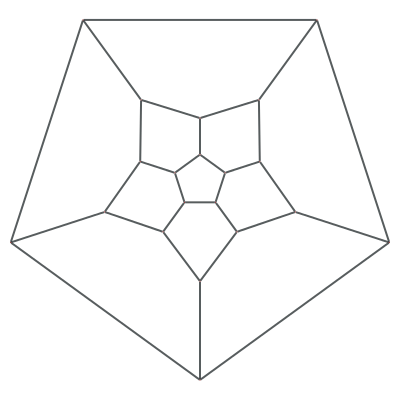}&
2
\includegraphics[height=3cm, origin=c]{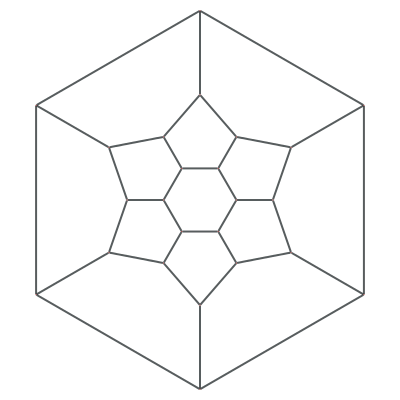}&
3
\includegraphics[height=3cm, origin=c]{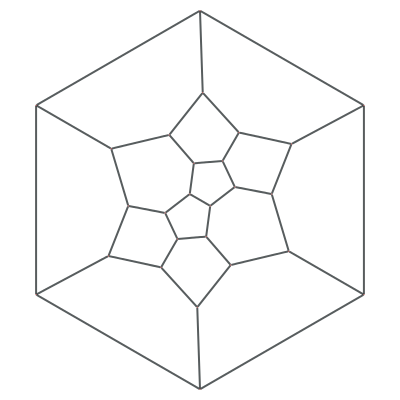}&
4
\includegraphics[height=3cm, origin=c]{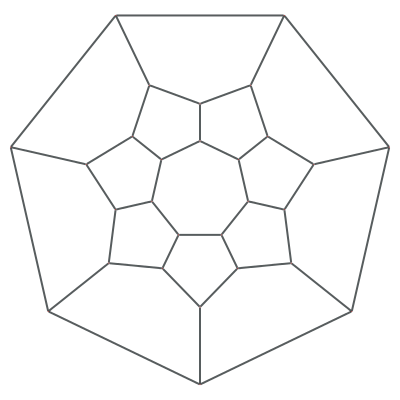}\\
\hline
5
\includegraphics[height=3cm, origin=c]{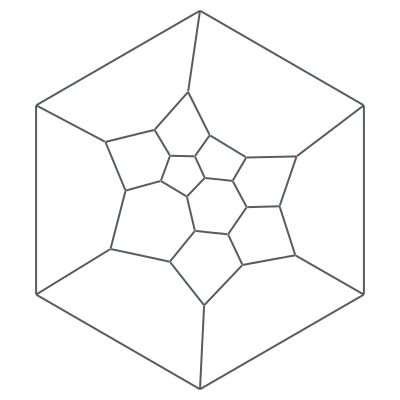}&
6
\includegraphics[height=3cm, origin=c]{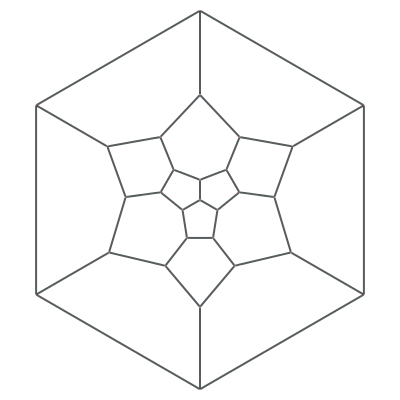}&
7
\includegraphics[height=3cm, origin=c]{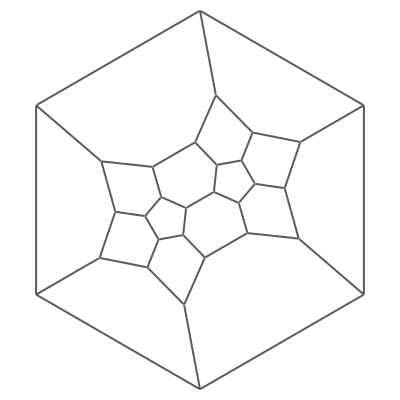}&
8
\includegraphics[height=3cm, origin=c]{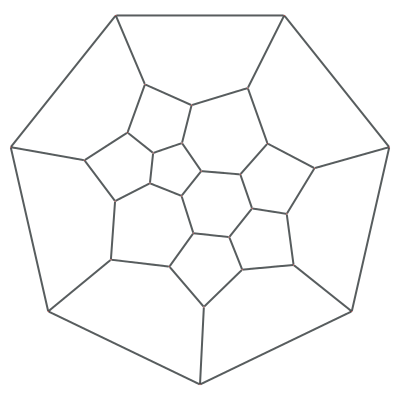}\\
\hline
9
\includegraphics[height=3cm, origin=c]{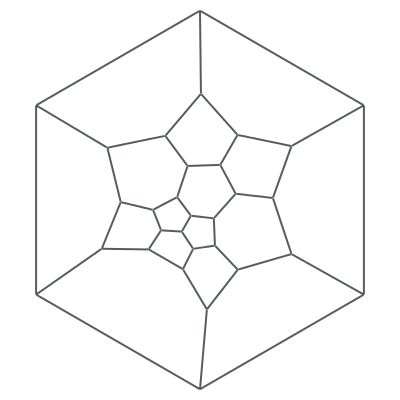}&
10
\includegraphics[height=3cm, origin=c]{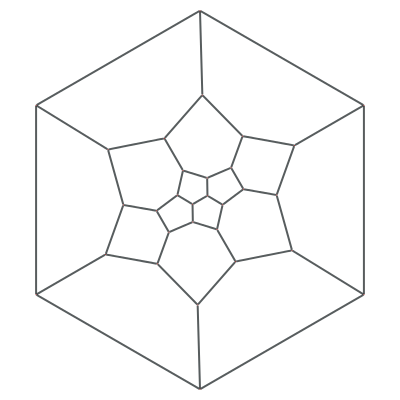}&
11
\includegraphics[height=3cm, origin=c]{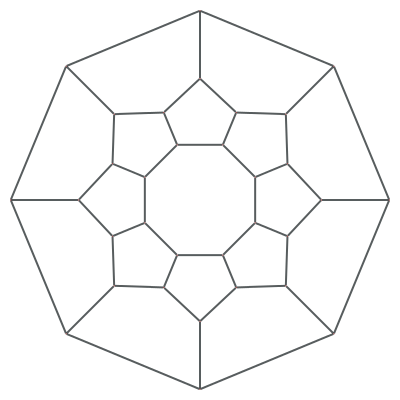}&
12
\includegraphics[height=3cm, origin=c]{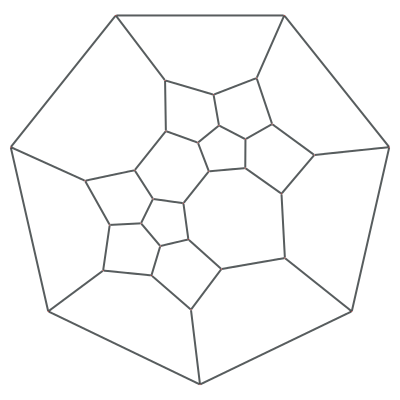}\\
\hline
13
\includegraphics[height=3cm, origin=c]{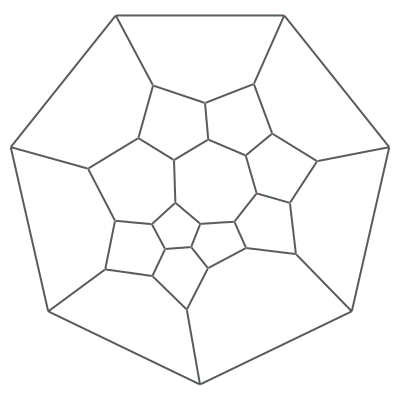}&
14
\includegraphics[height=3cm, origin=c]{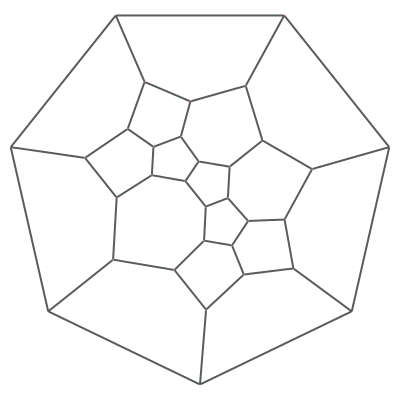}&
15
\includegraphics[height=3cm, origin=c]{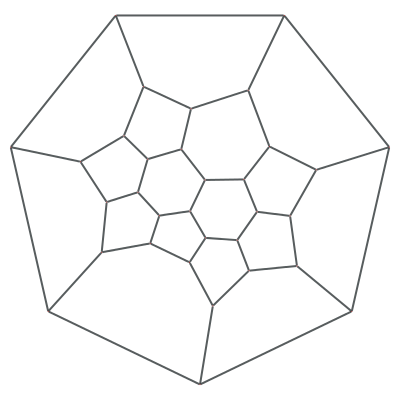}&
16
\includegraphics[height=3cm, origin=c]{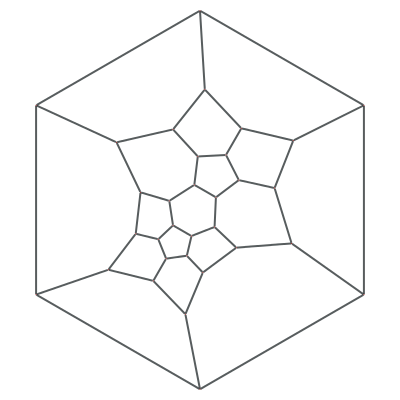}\\
\hline
17
\includegraphics[height=3cm, origin=c]{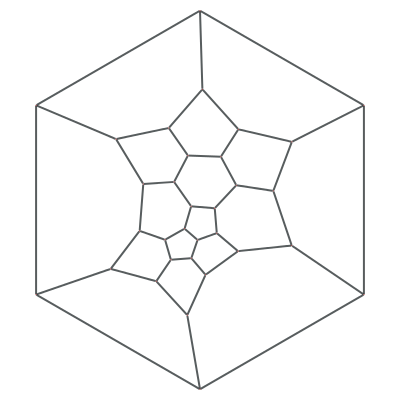}&
18
\includegraphics[height=3cm, origin=c]{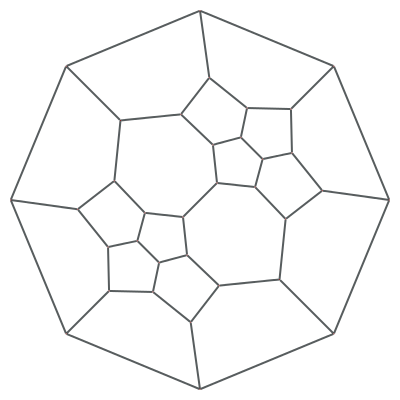}&
19
\includegraphics[height=3cm, origin=c]{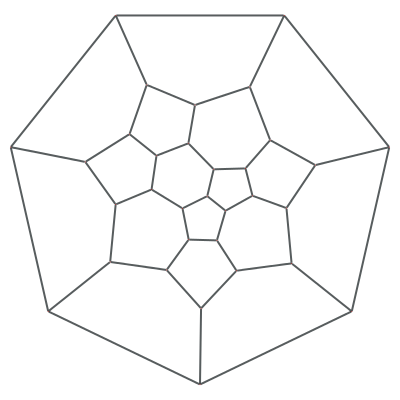}&
20
\includegraphics[height=3cm, origin=c]{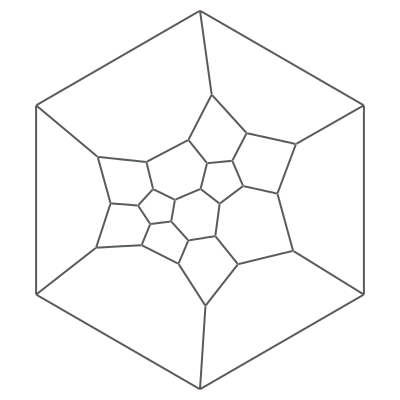}\\
\hline
\end{tabular} 
\newpage 
\begin{tabular}{|m{0.22\textwidth} | m{0.22\textwidth} | m{0.22\textwidth} |  m{0.22\textwidth} | } 
 \hline 
21
\includegraphics[height=3cm, origin=c]{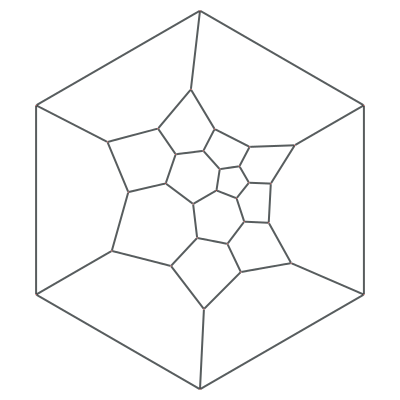}&
22
\includegraphics[height=3cm, origin=c]{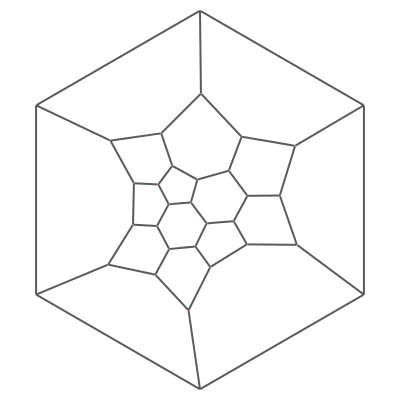}&
23
\includegraphics[height=3cm, origin=c]{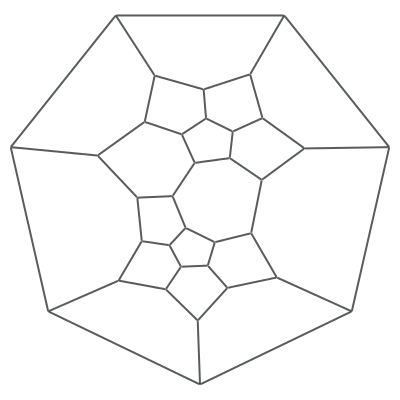}&
24
\includegraphics[height=3cm, origin=c]{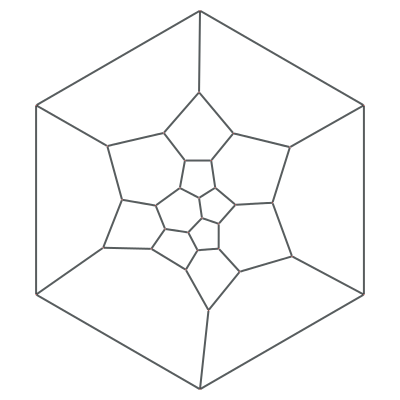}\\
\hline
25
\includegraphics[height=3cm, origin=c]{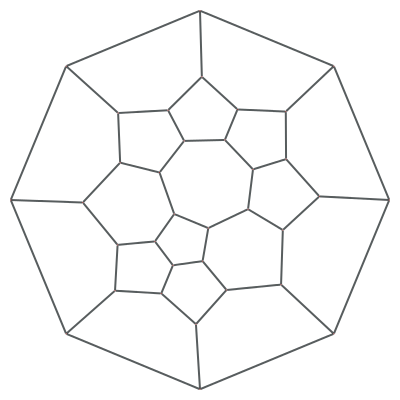}&
26
\includegraphics[height=3cm, origin=c]{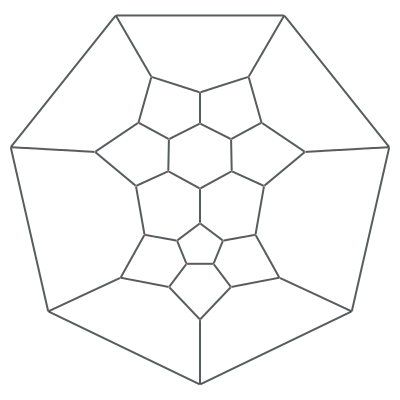}&
27
\includegraphics[height=3cm, origin=c]{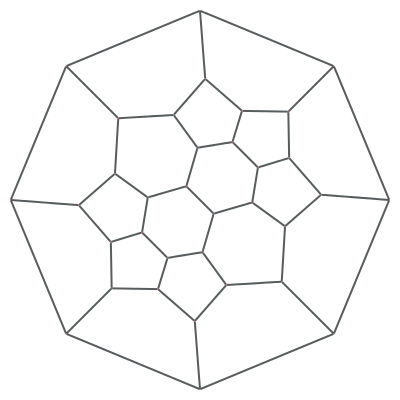}&
28
\includegraphics[height=3cm, origin=c]{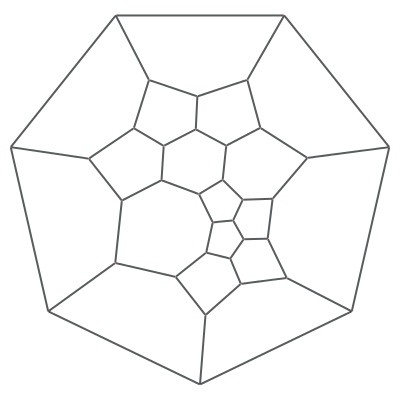}\\
\hline
29
\includegraphics[height=3cm, origin=c]{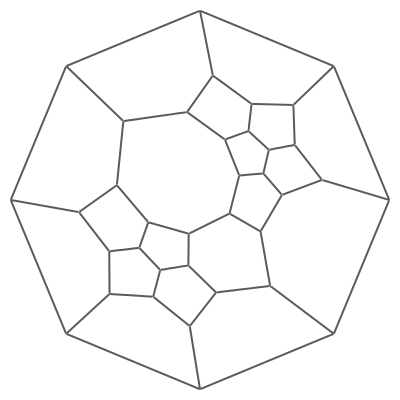}&
30
\includegraphics[height=3cm, origin=c]{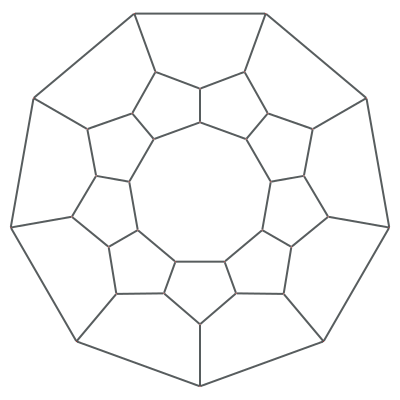}&
31
\includegraphics[height=3cm, origin=c]{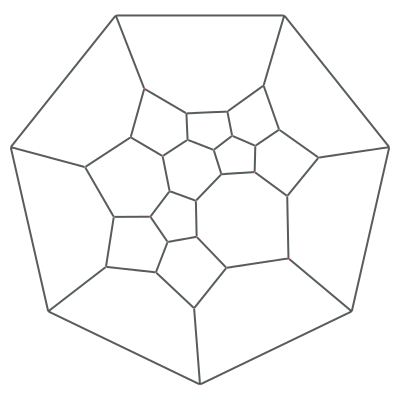}&
32
\includegraphics[height=3cm, origin=c]{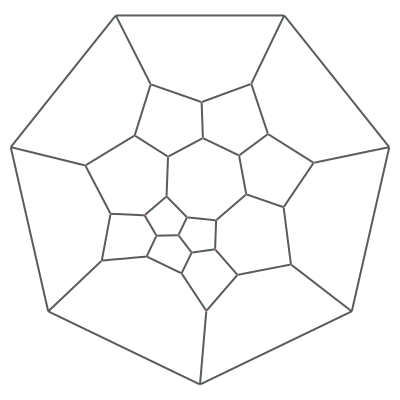}\\
\hline
33
\includegraphics[height=3cm, origin=c]{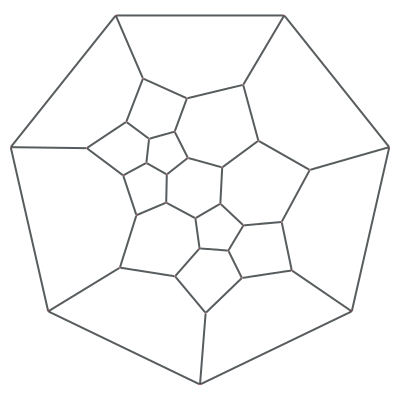}&
34
\includegraphics[height=3cm, origin=c]{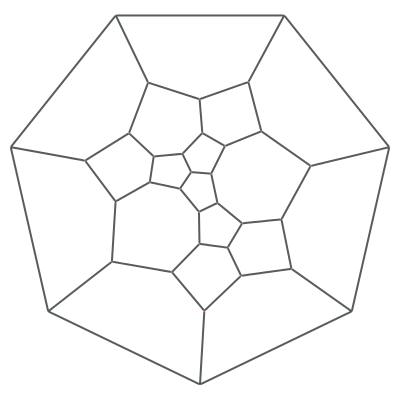}&
35
\includegraphics[height=3cm, origin=c]{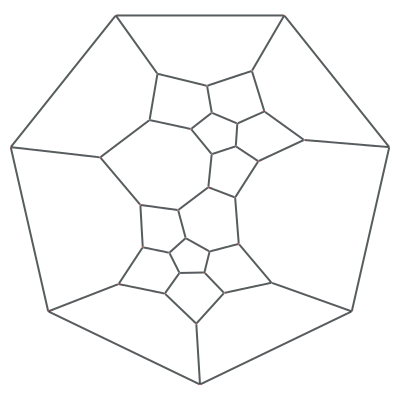}&
36
\includegraphics[height=3cm, origin=c]{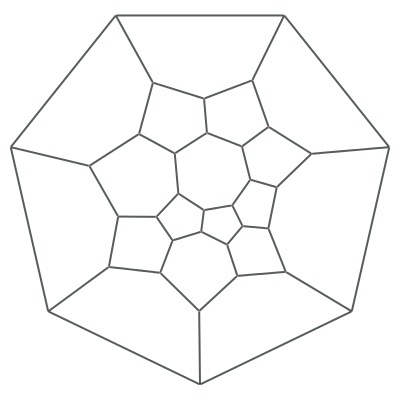}\\
\hline
37
\includegraphics[height=3cm, origin=c]{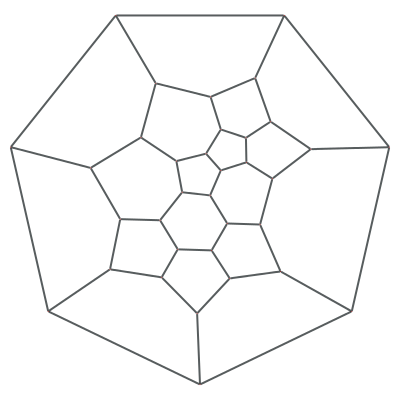}&
38
\includegraphics[height=3cm, origin=c]{38.png}&
39
\includegraphics[height=3cm, origin=c]{39.png}&
40
\includegraphics[height=3cm, origin=c]{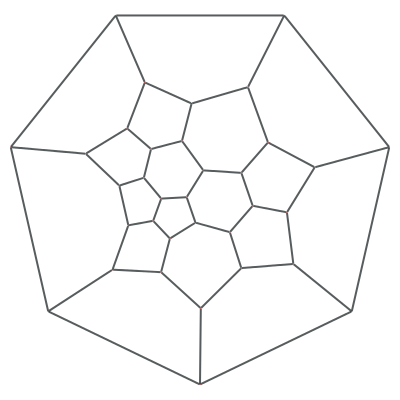}\\
\hline
\end{tabular} 
\newpage 
\begin{tabular}{|m{0.22\textwidth} | m{0.22\textwidth} | m{0.22\textwidth} |  m{0.22\textwidth} | } 
 \hline 
41
\includegraphics[height=3cm, origin=c]{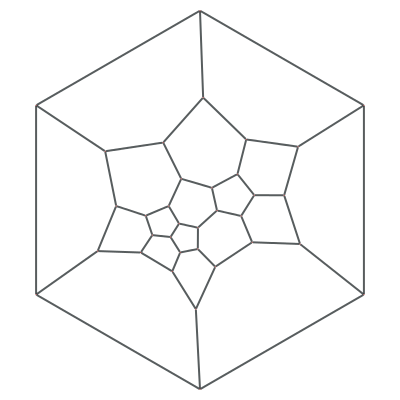}&
42
\includegraphics[height=3cm, origin=c]{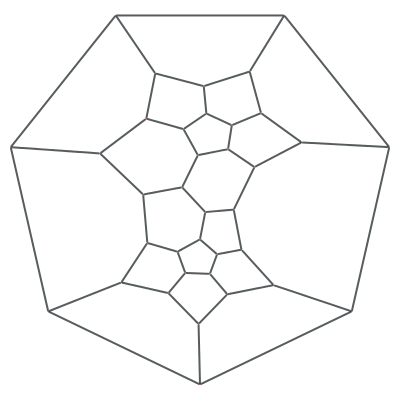}&
43
\includegraphics[height=3cm, origin=c]{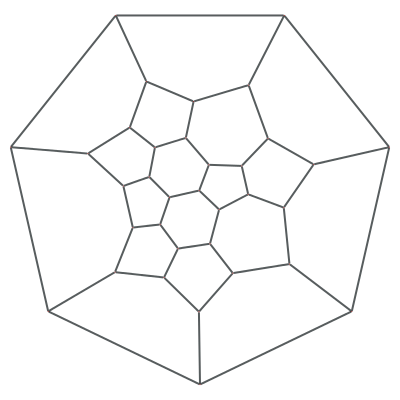}&
44
\includegraphics[height=3cm, origin=c]{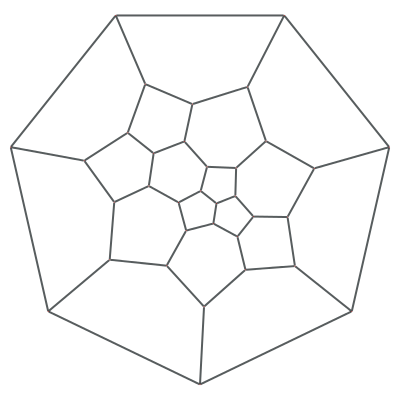}\\
\hline
45
\includegraphics[height=3cm, origin=c]{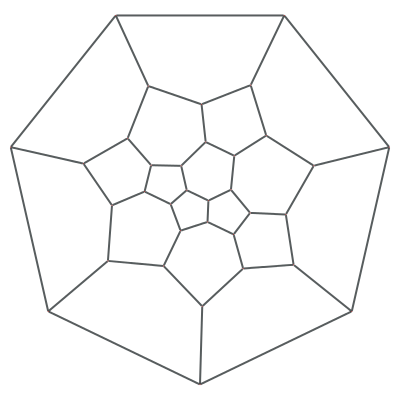}&
46
\includegraphics[height=3cm, origin=c]{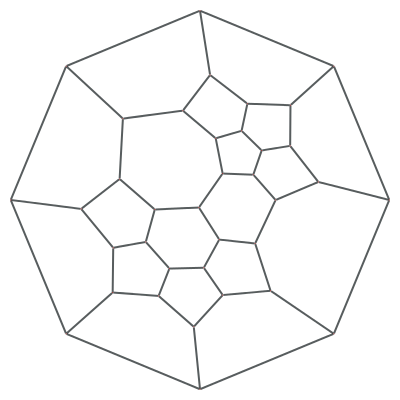}&
47
\includegraphics[height=3cm, origin=c]{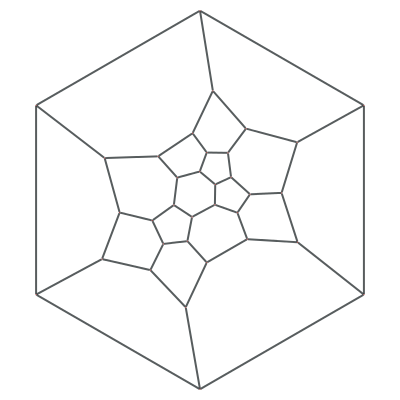}&
48
\includegraphics[height=3cm, origin=c]{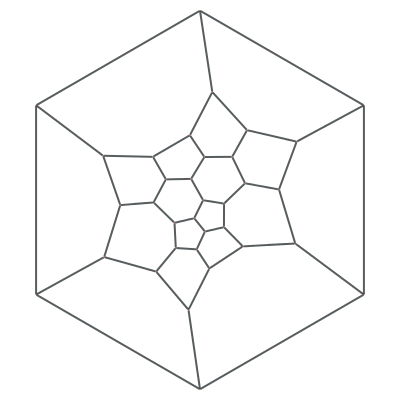}\\
\hline
49
\includegraphics[height=3cm, origin=c]{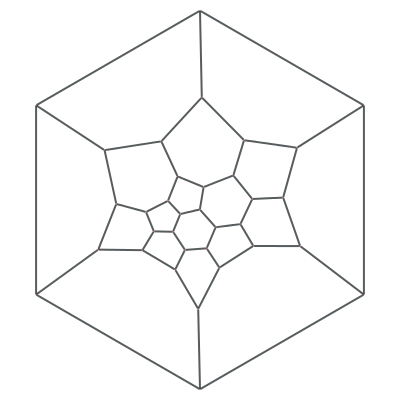}&
50
\includegraphics[height=3cm, origin=c]{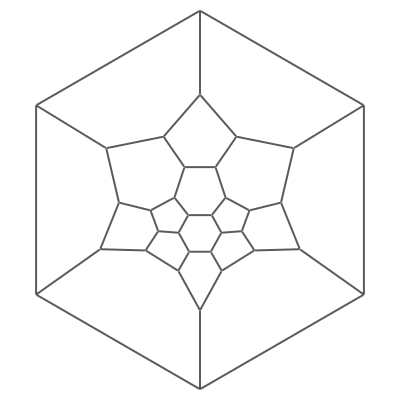}&
51
\includegraphics[height=3cm, origin=c]{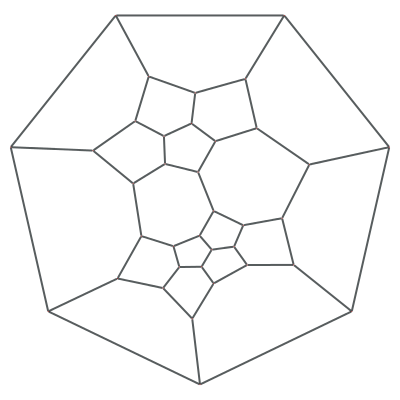}&
52
\includegraphics[height=3cm, origin=c]{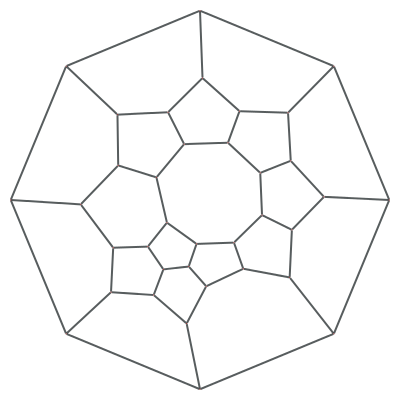}\\
\hline
53
\includegraphics[height=3cm, origin=c]{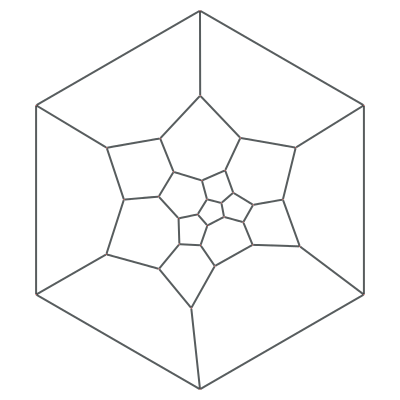}&
54
\includegraphics[height=3cm, origin=c]{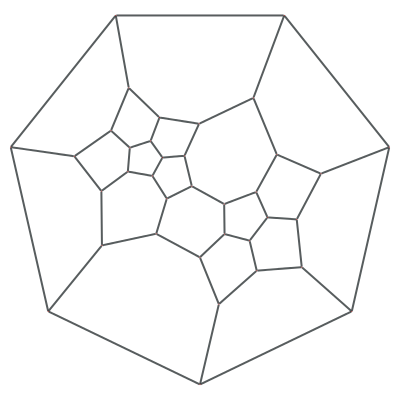}&
55
\includegraphics[height=3cm, origin=c]{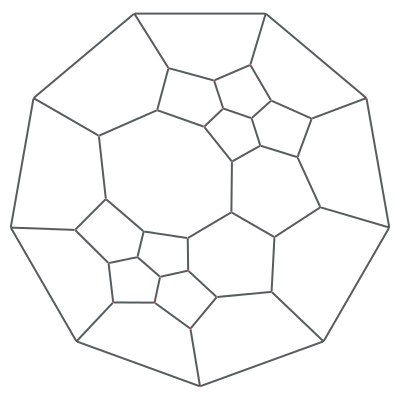}&
56
\includegraphics[height=3cm, origin=c]{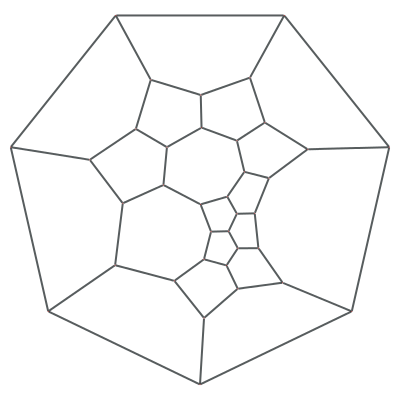}\\
\hline
57
\includegraphics[height=3cm, origin=c]{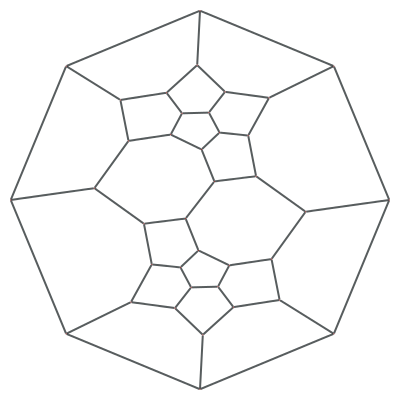}&
58
\includegraphics[height=3cm, origin=c]{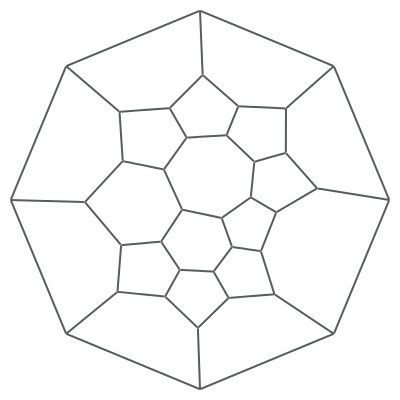}&
59
\includegraphics[height=3cm, origin=c]{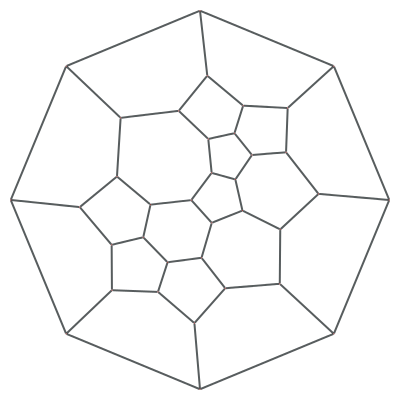}&
60
\includegraphics[height=3cm, origin=c]{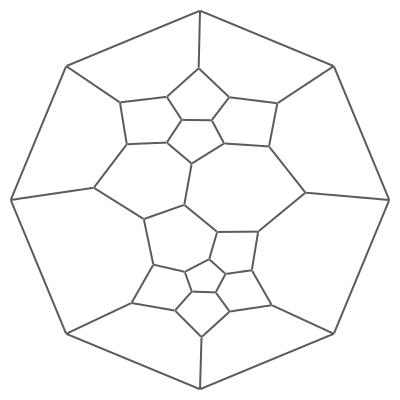}\\
\hline
\end{tabular} 
\newpage 
\begin{tabular}{|m{0.22\textwidth} | m{0.22\textwidth} | m{0.22\textwidth} |  m{0.22\textwidth} | } 
 \hline 
61
\includegraphics[height=3cm, origin=c]{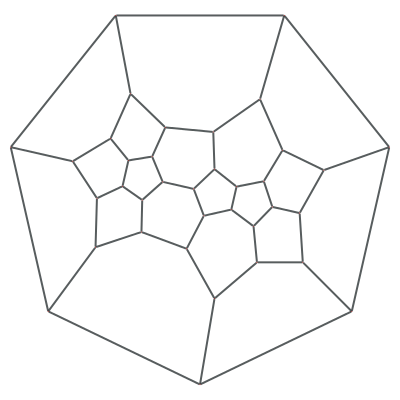}&
62
\includegraphics[height=3cm, origin=c]{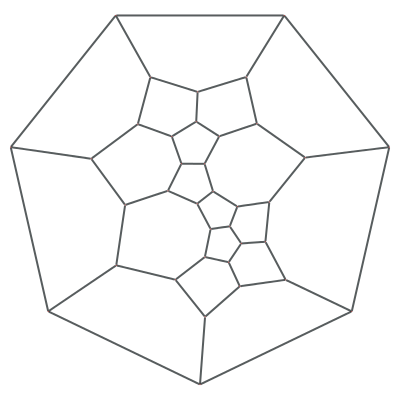}&
63
\includegraphics[height=3cm, origin=c]{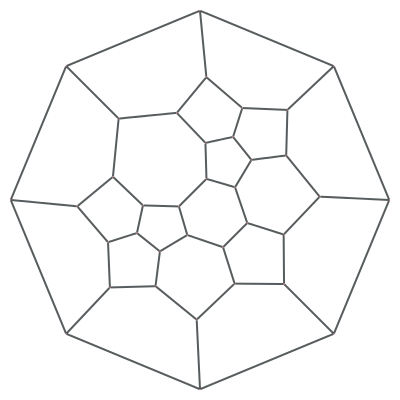}&
64
\includegraphics[height=3cm, origin=c]{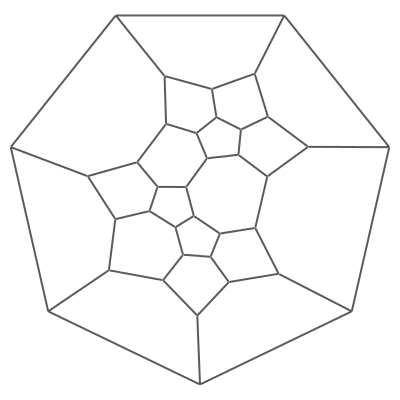}\\
\hline
65
\includegraphics[height=3cm, origin=c]{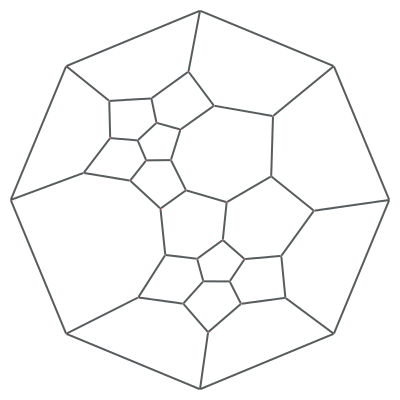}&
66
\includegraphics[height=3cm, origin=c]{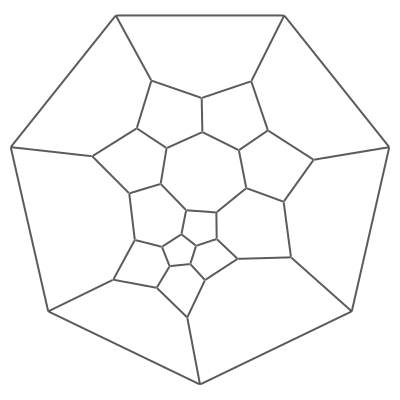}&
67
\includegraphics[height=3cm, origin=c]{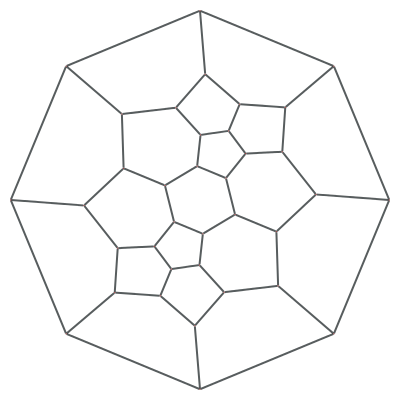}&
68
\includegraphics[height=3cm, origin=c]{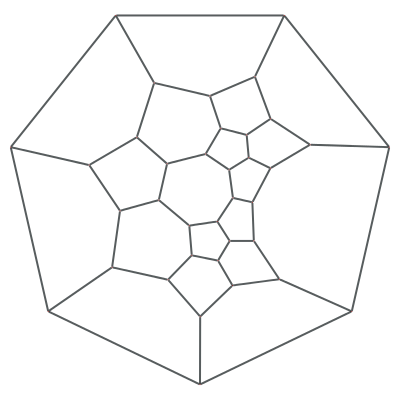}\\
\hline
69
\includegraphics[height=3cm, origin=c]{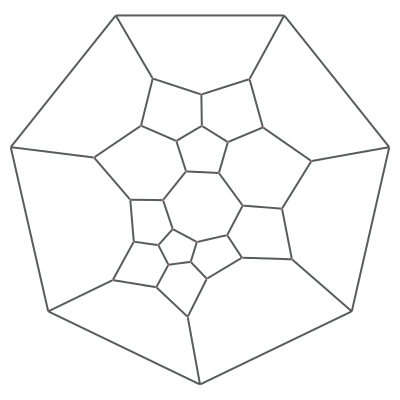}&
70
\includegraphics[height=3cm, origin=c]{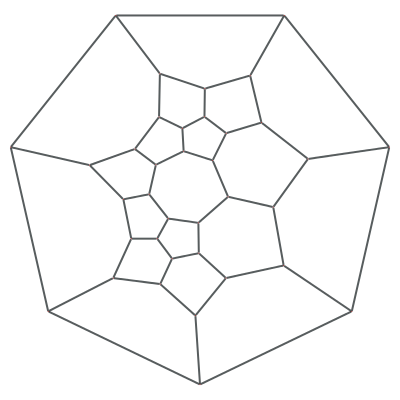}&
71
\includegraphics[height=3cm, origin=c]{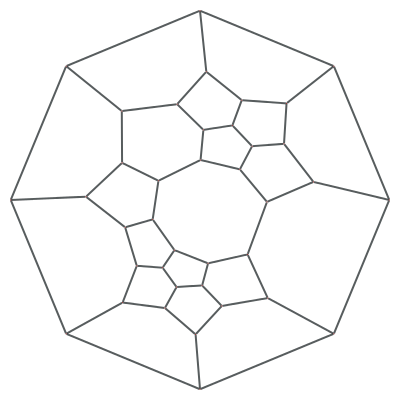}&
72
\includegraphics[height=3cm, origin=c]{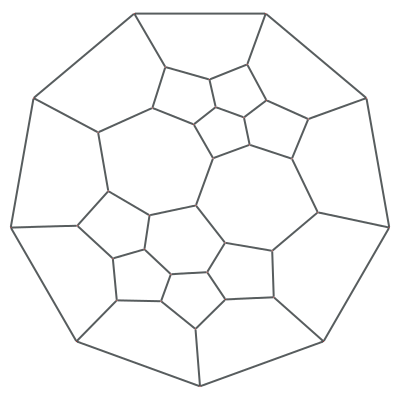}\\
\hline
73
\includegraphics[height=3cm, origin=c]{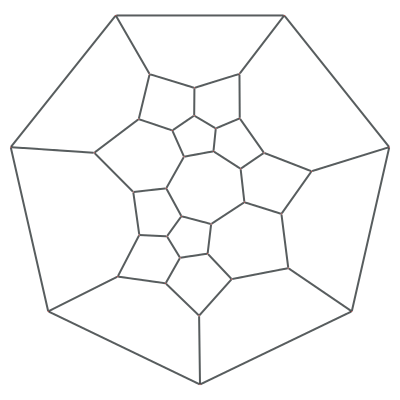}&
74
\includegraphics[height=3cm, origin=c]{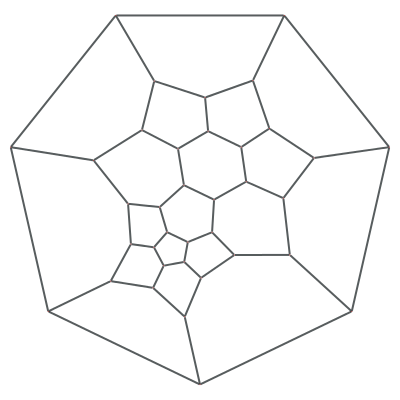}&
75
\includegraphics[height=3cm, origin=c]{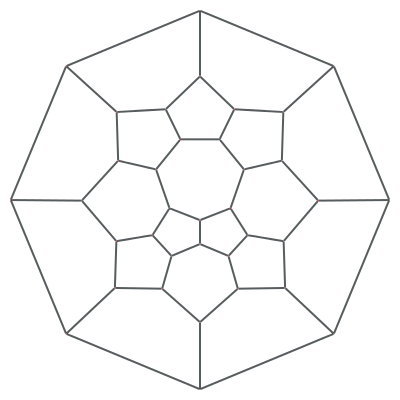}&
76
\includegraphics[height=3cm, origin=c]{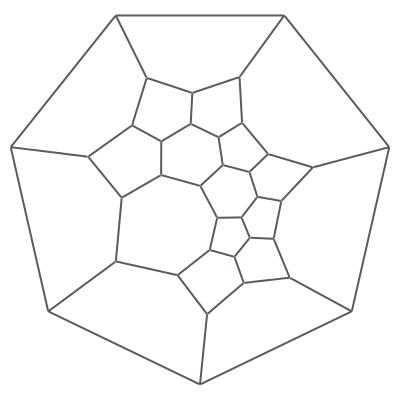}\\
\hline
77
\includegraphics[height=3cm, origin=c]{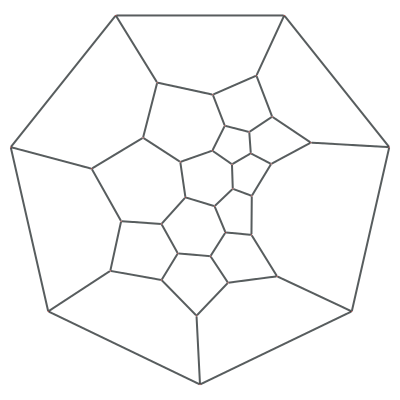}&
78
\includegraphics[height=3cm, origin=c]{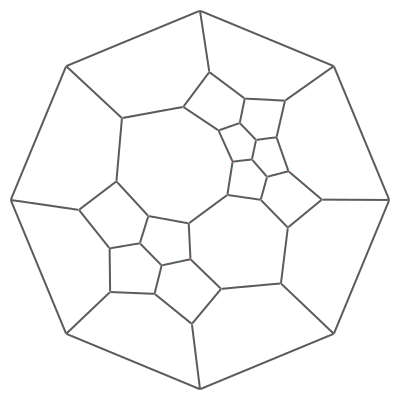}&
79
\includegraphics[height=3cm, origin=c]{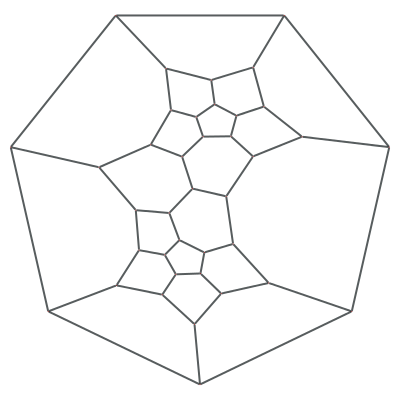}&
80
\includegraphics[height=3cm, origin=c]{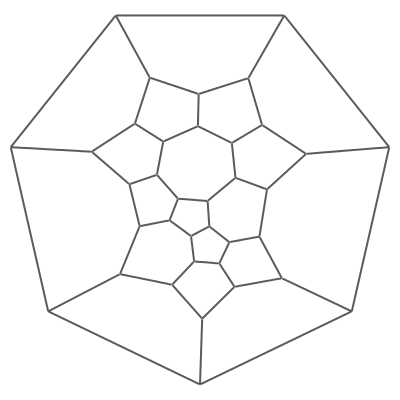}\\
\hline
\end{tabular} 
\newpage 
\begin{tabular}{|m{0.22\textwidth} | m{0.22\textwidth} | m{0.22\textwidth} |  m{0.22\textwidth} | } 
 \hline 
81
\includegraphics[height=3cm, origin=c]{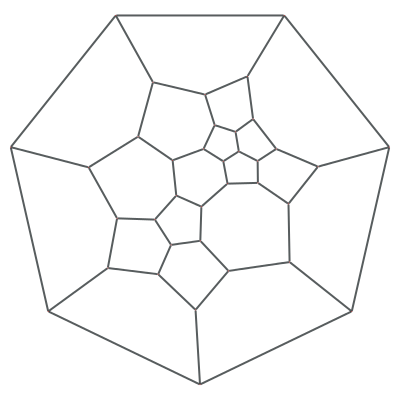}&
82
\includegraphics[height=3cm, origin=c]{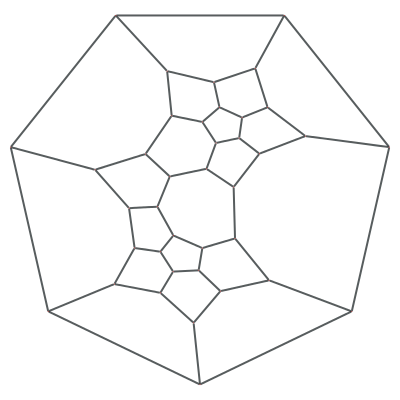}&
83
\includegraphics[height=3cm, origin=c]{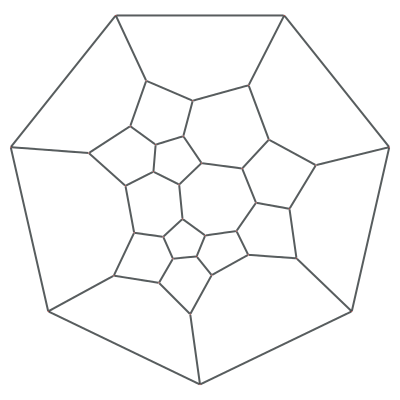}&
84
\includegraphics[height=3cm, origin=c]{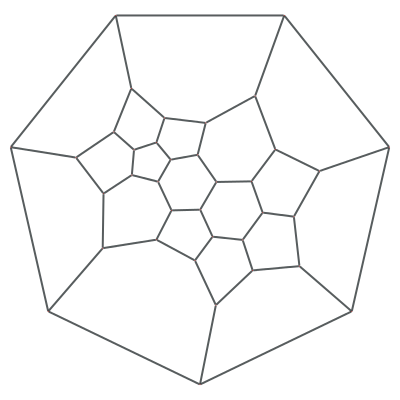}\\
\hline
85
\includegraphics[height=3cm, origin=c]{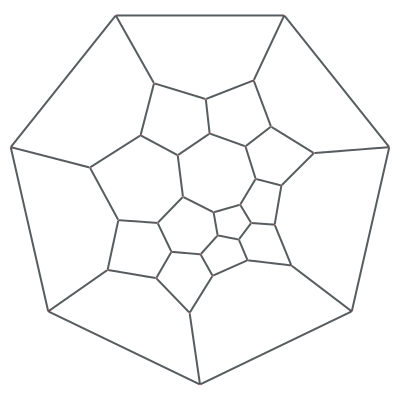}&
86
\includegraphics[height=3cm, origin=c]{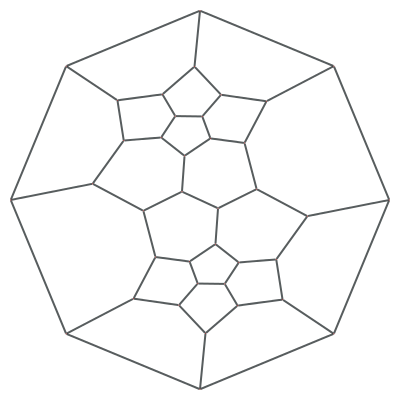}&
87
\includegraphics[height=3cm, origin=c]{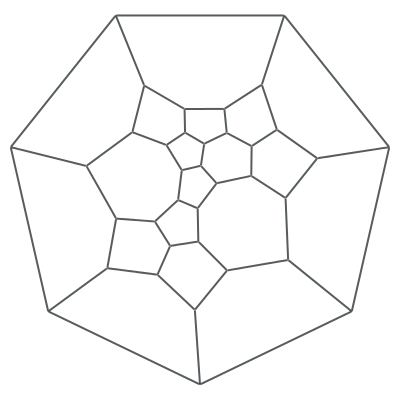}&
88
\includegraphics[height=3cm, origin=c]{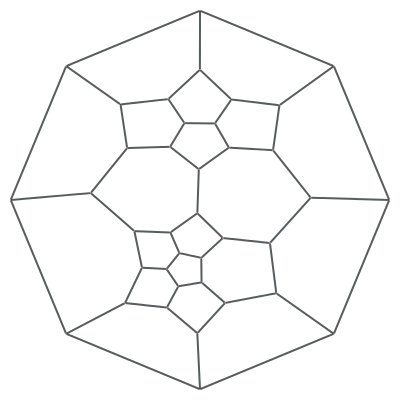}\\
\hline
89
\includegraphics[height=3cm, origin=c]{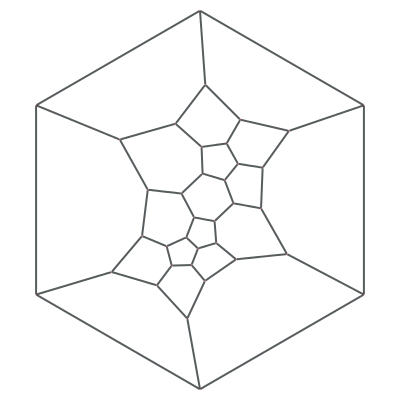}&
90
\includegraphics[height=3cm, origin=c]{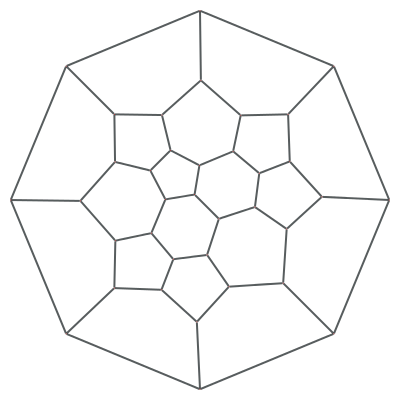}&
91
\includegraphics[height=3cm, origin=c]{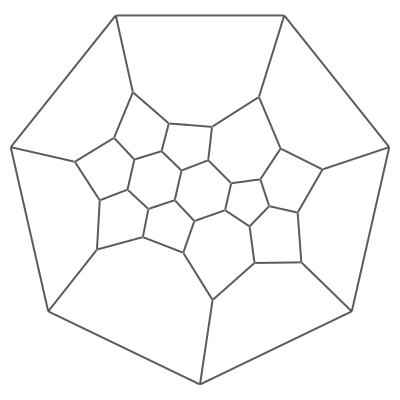}&
92
\includegraphics[height=3cm, origin=c]{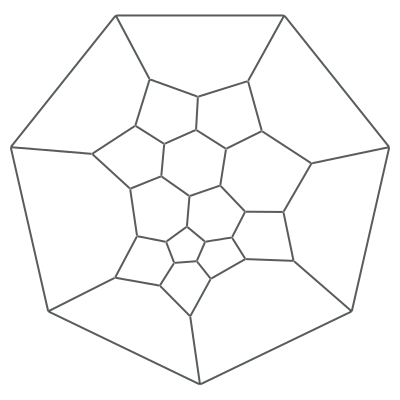}\\
\hline
93
\includegraphics[height=3cm, origin=c]{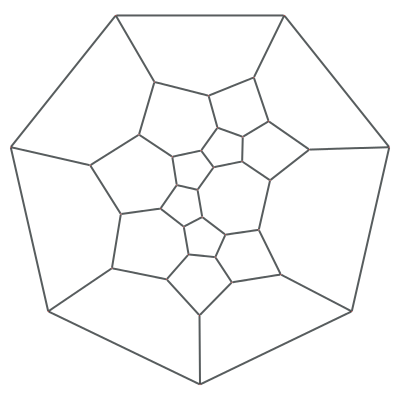}&
94
\includegraphics[height=3cm, origin=c]{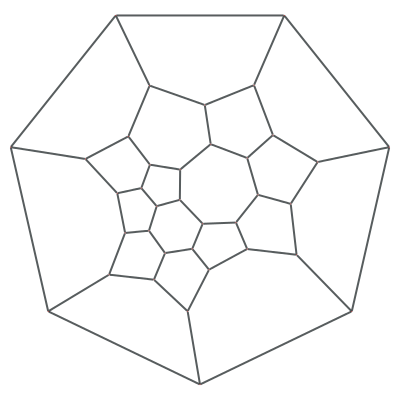}&
95
\includegraphics[height=3cm, origin=c]{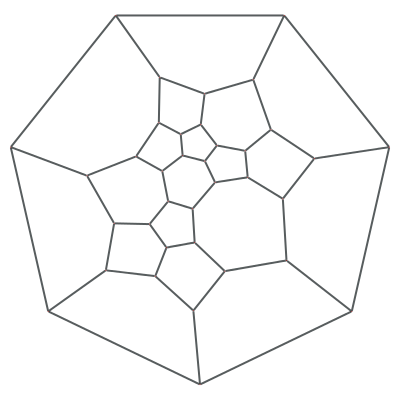}&
96
\includegraphics[height=3cm, origin=c]{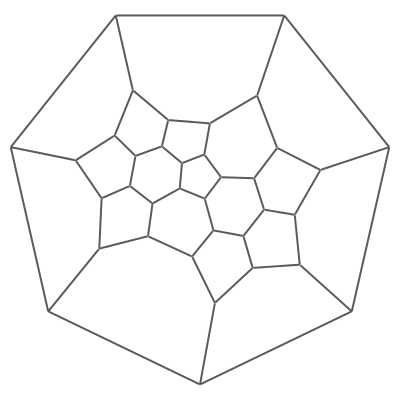}\\
\hline
97
\includegraphics[height=3cm, origin=c]{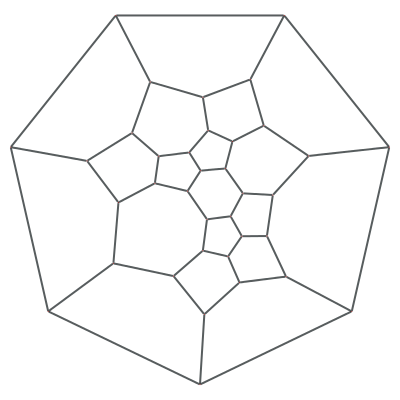}&
98
\includegraphics[height=3cm, origin=c]{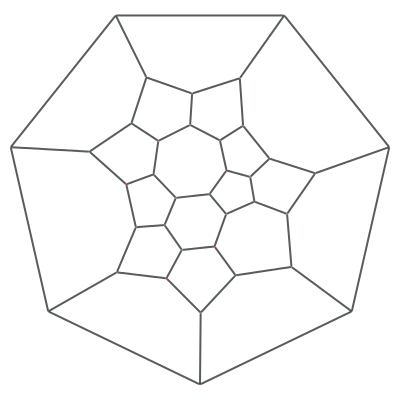}&
99
\includegraphics[height=3cm, origin=c]{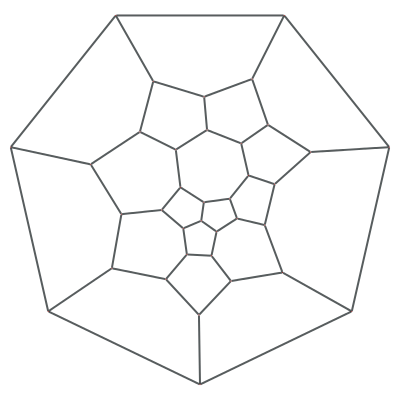}&
100
\includegraphics[height=3cm, origin=c]{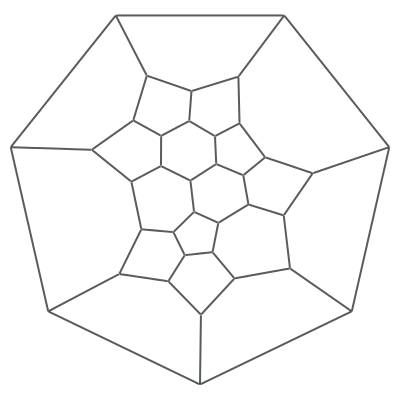}\\
\hline
\end{tabular}  
\newpage

\end{document}